\documentclass[12pt]{article}

\usepackage{amssymb}

\makeatletter \@addtoreset{equation}{section} \makeatother

\newcommand{\beq}{\begin{equation}}
\newcommand{\eeq}{\end{equation}}
\newcommand{\bea}{\begin{eqnarray}}
\newcommand{\eea}{\end{eqnarray}}

\newcommand{\req}[1]{(\ref{#1})}

\newtheorem{cor}{Corollary}[section]
\newtheorem{prop}{Proposition}[section]

\begin{document}

\title{Inverse Approach
 In The Study Of Ordinary Differential Equations.}

\title{
{\textbf{ Inverse Approach
 In The Study Of Ordinary Differential Equations}
} \\
\vspace{1ex} }

\date{}

\author{Rafael Ramirez$^{1}$ and
Natalia Sadovskaia$^2$
\\
\vspace{0.001ex} \\
\normalsize  $^1$ Departament d'Enginyeria Inform\`{a}tica i Matem\`{a}tiques, \\
\normalsize  Universitat Rovira i Virgili, \\
\normalsize  Avinguda dels Pa\"{\i}sos Catalans 26, \\
\normalsize  43007 Tarragona, Spain \\
\normalsize  \\
\normalsize  $^2$Departamento de Matem\'atica Aplicada II, \\
\normalsize  Universitat Polit\`ecnica de Catalunya \\
\normalsize  C. Pau Gargallo,5 \\
\normalsize  08028, Spain \\
}

\maketitle

\begin{abstract}

 We extend   the Eruguin result exposed in the paper "Construction of the whole set of ordinary differential equations
with a given integral curve" published in 1952 and construct a
differential system in $\Bbb{R}^N$ which admits a given set of the
partial integrals, in particular we study the case when theses
functions are polynomials. We construct a non-Darboux integrable
planar polynomial system of degree $n$ with one invariant
irreducible algebraic curve $g(x,y)=0$. For this system we analyze
the Darboux integrability, Poincare's problem and 16th's Hilbert
problem for algebraic limit cycles.

 {\textbf{ Mathematics Subject Classification (2000),}} 14P25,
34C05, 34A34.

\vskip0.25cm

\end{abstract}

\section{Introduction.}

Nonlinear ordinary differential equations appear in many branches
of applied mathematics, physics and, in general, in applied
sciences.

 By definition a real autonomous
 differential system  is a differential system of the form
$$\dot{\textbf{x}}=\textbf{v}(x),\quad \textbf{x}\in\mathbb{R}^N$$
where the dependent variables $x=(x^1,x^2,...,x^N)$ are real , the
independent variable (time $t$) is real and functions
$\bold{v}(x)=(v^1(x),...,v^N(x))$  are continuous functions in
$\textsc{D}\subset{\mathbb{R}^N}.$

{\textbf{ Definition 1.1}} {\it The smooth function $g$ and the
relation $g(x)=0$ are said partial integral and invariant relation
of the vector field $\bold{v}(x)$ respectively  if
$$dg(\bold{v})|_{g(x)=0}=0 .$$}

  In this paper we are mainly interested in to study the differential
system which possess a given set of invariant  relations.

It is always helpful to look at this problem from another point of
view. In this paper, we take an alternative viewpoint of starting
with a given set of   invariant  relations and determining the
form of the system which has such a set as invariant set.
[\cite{Gal1},\,\cite{Gal2},\cite{Chr},\,\cite{Llib2},\,
\cite{Koo},\,\cite{Baut},\,\cite{Dol},\,\cite{Darb},\,\cite{Jou},\,\cite{Sad1},\,\cite{Sad2},\,
\cite{Ram1},\, \cite{Ram2}].

This approach was first developed by Eruguin in the paper
"Construction of the whole set of ordinary differential equations
with a given integral curve" published in 1952 \cite{Eru}. In that
article the author stated and solved the problem of constructing a
planar vector field for which the given curve is its invariant. It
is important to observe that Eruguin considered only one curve,
moreover he didn't require that this curve was necessarily
algebraic.

Eruguin proved that the most general planar vector field $\bold{v}$ for which the given curve

\begin{equation}\label{11}
  g(x,y)=0
\end{equation}
is its invariant curve generates the following differential equations

\begin{equation}\label{12}
\left\{%
\begin{array}{cc}
\dot{x}=\nu (x,y)\{g,\,x\}+a(x,y)\\
\dot{y}=\nu (x,y)\{g,\,y\}+b(x,y)
\end{array}%
 \right.\end{equation}
where and $\nu ,\,a,\,b$ are  functions which we determine from
the condition:

\begin{equation}\label{13}
 dg(\bold{v})=\Phi (x,y),\quad \Phi|_{g=0}=0
\end{equation}
where
$$\{g,\,f\}\equiv{\partial_xg\partial_yf-\partial_yg\partial_xf}$$

 These Eruguin ideas were applied in different
areas.  In particular Zubov in \cite{ Zub} constructed the planar
system with a given region of stability.

Zubov constructed the following vector field
\[
\left\{%
\begin{array}{cc}
\dot{x}=f\gamma \{g_1,\,x\}+g_1(\gamma\{f,\,x\}+g_1\varphi d_1)\\
\dot{y}=f\gamma \{g_1,\,y\}+g_1(\gamma\{f,\,y\}+g_1\varphi
d_2)
\end{array}%
 \right.\]
 where $f,\,\gamma
,\,d_1,\,d_2,\,\varphi $ are arbitrary functions which he choose
in such a away that
$$d_1(g_1\{f,\,y\}-f\{g_1,\,y\})+d_2(g_1\{f,\,x\}-f\{g_1,\,x\})=1$$
Under this condition Zubov proved that the following relations
holds
\[
\left\{%
\begin{array}{cc}
dg_1(\bold{v})=g_1(\gamma\{f,\,g_1\})+g^2_1\varphi (d_1\{g_1,\,y\}-d_2\{g_1,\,x\})\\
dG_1(\bold{v})=\varphi G_1, \quad{\ln
G_1\equiv{{\displaystyle\frac{h}{g_1}}}}
\end{array}%
 \right.\]

Galliulin in \cite{Gal2} determines the most general vector field
in $\mathbb{R}^N$ for which the given relations
$$g_j(t,x^1,\,x^2,...,x^N)=0,\quad j=1,2,..,S<\,N$$
are the invariant relations, where $g_1,g_2,..,g_S$ are smooth
independent functions. The constructed system is the following

$$\dot{\textbf{x}}=\displaystyle\frac{1}{\Gamma}\sum_{i,j=1}^S\Gamma_{ij}(\Phi_j-
\partial_{t}g_j)\,\textbf{grad}\,{g_i}+\textbf{Y}. $$
where $\textbf{x}=col(x^1,x^2,...,x^N),$  $\textbf{Y}$ is an
arbitrary vector orthogonal to the vectors
$$\textbf{grad}\,{g_j}=col(\partial_1g_j,\partial_2g_j,...,\partial_Ng_j),\quad j=1,
...,S,\quad \partial_k\equiv{\frac{\partial}{\partial{x^k}}},
$$ $\Gamma$ is the Grama determinant,  $\Gamma_{ij}$ are
the minors of $\Gamma$ and $\Phi_1,...\Phi_S$
 are arbitrary functions:
 $$\Phi_j|_{g_j=0}=0,\quad j=1,2,..,S.$$

The aim of this paper is to extend the Eruguin-Galliulin ideas to
the case when the number of the given invariant relations  is
bigger than $N-1.$ The results which we expose have been
systematically developed in \cite{Sad1}.

\section{ Definitions and statement of the main results}

In this section we constructed the most general stationary
differential system from  the given set of partial integrals.

First of all we introduce the following concept and notations
which we shall use below.

{\textbf{ Definition 2.1}}\,

{\it We call the vector field:

\begin{equation}\label{21}
\textbf{v}=
-\displaystyle\frac{1}{\Upsilon}\left|\begin{array}{rrrrr}
dg_1(\partial_{1})&dg_1(\partial_{2})&\ldots &dg_1(\partial_{N})&\Phi_1 \\
dg_2(\partial_{1})&dg_2(\partial_{2})&\ldots &dg_2(\partial_{N})&\Phi_2 \\
\vdots           &\ldots                    &\vdots&\vdots\\
dg_M(\partial_{1})&dg_M(\partial_{2})&\ldots &dg_M(\partial_{N})&\Phi_M \\
dg_{M+1}(\partial_{1})&dg_{M+1}(\partial_{2})&\ldots &dg_{M+1}(\partial_{N})&\lambda_{M+1} \\
\vdots           &\ldots                    &\vdots&\vdots\\
dg_N(\partial_{1})&dg_N(\partial_{2})&\ldots &dg_N(\partial_{N})&\lambda_N \\
\partial_1&\partial_2&\ldots&\partial_N&0
\end{array}
\right|
\end{equation}
 the  Eruguin-Galliulin Vector Fields,}
where $g_1,\,g_2,\,...,g_{N}$ are smooth functions,
$\Phi_1,\Phi_2,...,\Phi_M$ are the Eruguin functions and
$\lambda_{M+1},\lambda_{M+2},...,\lambda_{N}$ are arbitrary
functions:

\begin{equation}\label{22}
\left\{%
 \begin{array}{cc}
 dg_k(\bold{v})=\Phi_k,\quad \Phi_k|_{g_k=0}=0,\quad k=1,..,M,\\
\\
dg_j(\bold{v})=\lambda_j,\quad j=M+1,..,N,
\end{array}%
\right.
\end{equation}

$dg_1,dg_2,...,dg_M$ are given independents 1-forms and
$dg_{M+1},dg_{M+2},...,dg_N$ are arbitrary 1-forms which we choose
in such a way that

\begin{equation}\label{23}
\Upsilon\equiv\left|\begin{array}{rrrr}
dg_1(\partial_{1})&dg_1(\partial_{2})&\ldots   &dg_1(\partial_{N})\\
dg_2(\partial_{1})&dg_2(\partial_{2})&\ldots   &dg_2(\partial_{N})\\
            \vdots&            \ldots&\ldots   &\vdots\\
 dg_N(\partial_{1})&dg_N(\partial_{2})&\ldots
&dg_N(\partial_{N})
\end{array}
\right| \equiv{\{g_1,\,g_2,\,...,g_{N}\}}\ne{0}
\end{equation}

{\it The functions $\Phi_1,\,\Phi_2,...,\Phi_M$ we call the
Eruguin functions}\,\cite{Gal2}.

We can identify the vector field \req{21} with the first order
differential system

\begin{equation}\label{24}
 \dot{\bold{x}}=\Upsilon M^{-1}\bold{w},
\end{equation}
 where $M$ and $\bold{w}$
are the matrices:
\[
\left\{%
 \begin{array}{cc}
 \textsc{M}=\Big((dg_j(\partial_k))_{j,k=1,2,..,N}\Big),\\
 \textbf{w}=col(\Phi_1,\,...,\Phi_M,\,\lambda_{M+1},...,\lambda_N).
\end{array}
\right.
\]

It is easy to show that the system \req{24} admits the equivalent
representation
\begin{equation}\label{25}
\dot{x}^j= \Phi_1\{x^j,g_2,..g_{M+1}...,g_N\}
+\ldots+\Phi_M\{g_1,..,x^j,\,g_{M+1},..,g_N\}+ Y^j,
\end{equation}
 where
$$Y^j=\lambda_{M+1}\{g_1..,g_{M},\,x^j,\,g_{M+2}..,g_N\,\}...+\lambda_{N}\{g_1..,g_M,g_{M+1}..,g_{N-1},\,x^j\},$$
 $j=1,2,.., N. $

Clearly, the vector $\textbf{Y}=col(Y^1,\,Y^2,...,Y^N)$ is
orthogonal to the vectors $\textbf{grad}{g_j},\,j=1,2,..,M,$ hence
we obtain the Galliulin result \cite{Gal2}.

 {\textbf{ Example 2.1}}\,

 We shall construct the Eruguin-Galliulin vector field
for the case when the arbitrary functions
$$g_{M+1},\,g_{M+2},\,...,g_{M+K},\quad N=M+K$$
are such that

\begin{equation}\label{26}
 \left\{%
 \begin{array}{cc}
 dg(\bold{v})=L\,g,\\
dg_{M+1}(\bold{v})=Lg_{M+1}+L_1g\\
dg_{M+2}(\bold{v})=Lg_{M+2}+L_1g_{M+1}+L_2g\\
\vdots\\
dg_{M+K}(\bold{v})=Lg_{M+K}+L_1g_{M+K-1}+...+L_Kg
\end{array}
\right.
\end{equation}

 where $L_1,\,L_2,...,L_K,\,L$ are arbitrary functions and
 $$g=\prod_{j=1}^Mg^{\tau_j}_j,\quad\tau_j\in{\textsc{C}}$$

By introducing the functions $G_1,\,G_2,\,...,G_{K}:$
$$g_{M+j}=G_jg,\quad j=1,2,..,K,$$
 we obtain
\[
\left\{%
 \begin{array}{cc}
dG_{1}(\bold{v})=L_1\\
dG_{2}(\bold{v})=L_1G_{1}+L_2\\
\vdots\\
dG_{K}(\bold{v})=L_1G_{K-1}+...+L_k
\end{array}
\right.
\]
 Clearly, the arbitrary functions
$\lambda_{M+1},\,\lambda_{M+2},...,\lambda_{N}=\lambda_{M+K}$ in
this case we determine as follow
\[
\left\{%
\begin{array}{cc}
\lambda_{M+1}=g(LG_{1}+L_1)\\
\lambda_{M+2}=g(LG_2+L_1G_{1}+L_2)\\
\vdots\\
\lambda_{M+K}=g(LG_K+L_1G_{K-1}+...+L_K)
\end{array}%
\right.
\]

 Let us introduce the
1-forms  $\omega_1,\,\omega_2,...,\omega_K:$
\[
\left\{%
\begin{array}{cc}
dG_1=\omega_1\\
dG_2=G_1\omega_1+\omega_2\\
\vdots\\
dG_{K}=G_{K-1}\omega_1+...+\omega_K.
\end{array}%
\right.
\]

 After some straightforward calculations we prove that
\[
\left\{%
\begin{array}{cc}
\omega_j=d\Upsilon_j,\\
\omega_j(\bold{v})=L_j,\quad j=1,2,..,K.
\end{array}%
\right.
\]
 Hence the functions $\Upsilon_1,\,\Upsilon_2,\,....\Upsilon_K$
are such that

\begin{equation}\label{27}
d\Theta=\Psi^{-1}dG,
\end{equation}

 where
\[
\left\{%
\begin{array}{cc}
 d\Theta=col(d\Upsilon_1,\,d\Upsilon_2,\,...,d\Upsilon_K),\\
dG=col(dG_1,\,dG_2,...dG_K),
\end{array}%
\right.
\]

\[
\Psi=\left(
\begin{array}{llllll}
 1 & 0 & 0 & 0 & \ldots & 0 \\
G_1  &1  &0   &0   &\ldots  &0\\
G_2 &G_1  &1   &0   &\ldots  &0\\
G_3 &G_2  &G_1   &1   &\ldots  &0\\
\vdots&\vdots &\vdots   &\vdots   &\ldots&0\\
G_K&G_{K-1}&G_{K-2}&\ldots& G_1 &1
\end{array}%
\right)\]

After the integration the system \req{27} we obtain

\[
\left\{%
\begin{array}{cc}
\Upsilon_1=G_1\\
\Upsilon_2=G_2-\frac{G^2_1}{2}\\
\Upsilon_3=G_3-G_1G_2+\frac{G^3_1}{3}\\
\Upsilon_4=G_4-G_1G_3+G^2_1G_2-\frac{G^4_1}{4}-\frac{G^2_2}{2}\\
\Upsilon_5=G_5-G_1G_4+G^2_1G_3-G^3_1G_2+\frac{G^5_1}{5!}+\frac{G^3_2}{3!}\\
\vdots
\end{array}%
\right.
\]

Hence, for the function $\Upsilon_j$  there are the equivalent
representations
$$\Upsilon_j=\sum_{m=1}^{K-1}\alpha_{jm}(x)g^{m-j},\quad j=1,...,K.$$
where $\alpha=(\alpha_{jm})$ is some matrix.

\begin{cor}

{\it Let us suppose that the functions

$$L,\,L_1,\,...,L_K$$ are such that

$$\sum_{j=0}^K\nu_jL_j=0,\quad L_0=L,$$
  then the constructed system \req{25},\req{26} admits the first integral}
$$F(x)=g^{\nu_0}\exp{\sum_{j=1}^K\nu_j\Upsilon_j}=\prod_{j=1}^Kg^{\nu_0\tau_j}_j\exp{\sum_{j=1}^K\nu_j\Upsilon_j},$$
where $\nu_0,\,\nu_1,...,\nu_K$ are  constants.
\end{cor}
 For the particular case when $N=2$ and the given curves

\begin{equation}\label{28}
  g_j(x)=0,\,j=1,2,..,M
\end{equation}
  are algebraic curve,  then from the Darboux's
 theory follows that
$$\Phi_j=K_j(x)\,g_j$$ thus , the condition on the existence the
first integral $F$ takes the form
$$\nu_0\sum_{j=1}^M\tau_jK_j+ \sum_{j=1}^K\nu_jL_j=0,$$

For the planar polynomial vector field this condition was deduced
in particular in \cite{Chr1}.  In this paper  the following
definition is given

\textbf{Definition of infinitesimal multiplicity}

{\it Let $f=0$ be an invariant algebraic curve of degree $n$ of a
polynomial vector field $X$ of degree $d$. We say that
$$F=f_0+f_1\epsilon+...+f_{k-1}\epsilon^{k-1}\in\textsc{C}[x,y,\epsilon
]/\epsilon^k$$ defines a generalized invariant algebraic curve of
order $k$ based on $f=0$ if $$f_0=f,\,f_1,\,...,f_{k-1}$$ are
polynomials in $\textsc{C}[x,y]$ of degree at most $n,$ and $F$
satisfies the equation
\begin{equation}\label{29}
X(F)=FL_F
\end{equation} for some polynomial
$$L_F=L_0+L_1\epsilon+...+L_{k-1}\epsilon^{k-1}\in\textsc{C}[x,y,\epsilon
]/\epsilon^k$$ which must necessarily be of degree at most $d-1$
in $x$ and $y.$ We call $L_F$ the cofactor of $F.$ Equations
\req{29} can be written as}
\[
\left\{%
\begin{array}{cc}
X(f_0)=f_0L_0\\
X(f_1)=f_1L_0+f_0L_1\\
\vdots \\
X(f_{k-1})=f_{k-1}L0+f_{k-2}L_1+...+f_0L_{k-1}
\end{array}%
\right.
\]
The vector field \req{25},\,\req{26} can be applied to extend the
concept of infinitesimal multiplicity for the polynomial  vector
field in $\mathbb{R}^N.$

 \begin{prop}

  {\it Let $$g_j(x)=0,\quad
 x=(x^1,x^2,...,x^N),\quad
j=1,2,..,M<\,N$$ are invariant relations  of a differential system
(S).

Assume that

$$\Upsilon=\{g_1,\,g_2,\,...,g_M,g_{M+1}...,g_{N}\}\ne{0},$$ for arbitrary smooth functions $g_{M+1},\,g_{M+2},..,g_N.$

Then the following statement hold:

  System ($\textsc{S}$) can be written as \req{25}.}
\end{prop}
\textbf{Proof.}

\, Suppose that

\begin{equation}\label{210}
\dot x=X(x)
\end{equation}

 is a differential system having  $g_1,\,g_2,...,g_M$ as partial integrals. Then taking

\[
\left\{%
\begin{array}{cc}
\Phi_j=\frac{1}{\Upsilon}d g_j( X),\quad j=1,2,..,M\\
\\
\lambda_{k}=\frac{1}{\Upsilon}d g_{k}(X),\quad k= M+1,M+2,..,N
\end{array}%
\right.
\]

 we get that the system \req{24}, or, what is the same,
\req{25} becomes system \req{210}. Note that in the definition
 of $\Phi_j$ and $\lambda_j$ we have used that
 $\{g_1,\,g_2,\,...,g_{N}\}\ne{0}.$

Now we shall study the case when the given number of partial
integrals is $S>N.$

 If $S=N$ then the differential system \req{25} takes the form

\begin{equation}\label{211}
\dot{x}^j=\Phi_1\{x^j...,g_M,...,g_N\}
+\ldots+\Phi_{N}\{g_1,...,g_M,..,x^j\}\quad j=1..,N.
\end{equation}

  \begin{prop}

 {\it The differential system \req{211} admits the complementary
invariant relation
$$ g_\nu(x)=0,\quad\nu=N+1,...,S$$
if and only if}
\end{prop}

\begin{equation}\label{212}
\left|\begin{array}{rrrrr}
dg_1(\partial_{1})&dg_1(\partial_{2})&\ldots &dg_1(\partial_{N})&\Phi_1 \\
dg_2(\partial_{1})&dg_2(\partial_{2})&\ldots &dg_2(\partial_{N})&\Phi_2 \\
\vdots           &\ldots                    &\vdots&\vdots\\
dg_M(\partial_{1})&dg_M(\partial_{2})&\ldots &dg_M(\partial_{N})&\Phi_M \\
dg_{M+1}(\partial_{1})&dg_{M+1}(\partial_{2})&\ldots &dg_{M+1}(\partial_{N})&\Phi_{M+1} \\
\vdots           &\ldots                    &\vdots&\vdots\\
dg_N(\partial_{1})&dg_N(\partial_{2})&\ldots &dg_N(\partial_{N})&\Phi_N \\
dg(\partial_{1})&dg(\partial_{2})&\ldots
&dg(\partial_{N})&\Phi_\nu
\end{array}%
\right|
\end{equation}
or, what is the same,

\begin{equation}\label{213}
\Phi_1\{g..,g_M..,g_N\}+..+\Phi_N\{g_1..,g_{M}..,g\}+\Phi_\nu\{g_1..,g_M..,g_N\}=0.
\end{equation}
We obtain the proof from the equality
$$dg_\nu(\bold{v})=\Phi_\nu ,$$
which in view of \req{21} coincides with \req{213}.

Below we shall use the following identity

\begin{equation}\label{214}
\begin{array}{cc}
\{f_1,f_2,...,f_{N-1},g_1\}\{g_2,g_3,...,g_N,G\}+\\
+\{f_1,f_2,...,f_{N-1},g_2\}\{g_1,g_3,...,g_N,G\}+...\\
+\{f_1,f_2,...,f_{N-1},g_N\}\{g_1,g_2,...,g_{N-1},G\}+\\
\{f_1,f_2,...,f_{N-1},G\}\{g_1,g_2,...,g_{N-1},g_N\} \equiv 0.
\end{array}%
\end{equation}

The proof follow by considering that \req{214} is equivalent to
the relation
\begin{equation}\label{215}
\left|\begin{array}{rrrrr} dg_1(\partial_{1})&dg_1(\partial_{2})&
\ldots   &dg_1(\partial_{N})&\{f_1,f_2,f_3,...,f_{N-1},g_1\}\\
dg_2(\partial_{1})& dg_2(\partial_{2})&\ldots
&dg_2(\partial_{N})&\{f_1,f_2,f_3,...,f_{N-
1},g_2\}\\
\vdots           &\ldots                    &\vdots&\vdots\\
dg_N(\partial_{1})&dg_N(\partial_{2})&\ldots
&dg_N(\partial_{N})&\{f_1,f_2,f_3,...,f_{N-1},g_N\}
 \\
dG(\partial_1)&dG(\partial_2)&\ldots&dG(\partial_N)&
\{f_1,f_2,f_3,...,f_{N-1},G\}
\end{array}%
\right|
 \equiv{0}
\end{equation}

 It is easy to show, in view of identity \req{214} that the Eruguin functions $\Phi_k,$ determined by
the formula
\begin{equation}\label{216}
\begin{array}{cc}
\Phi_k=\sum_{\alpha_1,\,\alpha_2,\,..,\alpha_{N-1}=1}^{S+N}
\{g_{\alpha_1},\,g_{\alpha_2}...,g_{\alpha_{N-1}},\,g_k\}\lambda_{{\alpha_1},{\alpha_2}...{\alpha_{N-1}}}(x),\\
k=1,2,..,S
\end{array}%
\end{equation}
 are the solutions
of \req{211}, where
$\lambda_{{\alpha_1},{\alpha_2}...{\alpha_{N-1}}},$ are  arbitrary
continuous  functions:
\begin{equation}\label{217}
 \Phi_k|_{g_k=0}=0,\quad k=1,2,....,S.
\end{equation}
 and $$g_{S+j}=x_j,\quad j=1,2,..,N$$
 The differential system \req{211} in this case takes
the form
\begin{equation}\label{218}
 \dot{x}^j=\sum_{\alpha_1,..,\alpha_{N-1}=1}^{S+N} \{
g_{\alpha_1},...,g_{\alpha_{N-1}},\,x^j\}\tilde{\lambda}_{{\alpha_1},{\alpha_2}...{\alpha_{N-1}}}(x)
\end{equation}

 where
$$\tilde{\lambda}_{{\alpha_1},{\alpha_2}...{\alpha_{N-1}}}(x)=\{g_1,g_2,g_3,...,g_N\}\lambda_{{\alpha_1},{\alpha_2}...{\alpha_{N-1}}}(x)
.$$

\begin{prop}

{\it Let $g_1(x),g_2(x),...,g_S(x)\quad S>N$ are partial integrals
of a differential system ($\textsc{S}$).

Assume that

$\Upsilon=\{g_1,\,g_2,...,g_{N}\}\ne{0},$

then the following statement hold:

  System ($\textsc{S}$) can be written as \req{218}.}

\end{prop}
\textbf{Proof.} In fact if we insert \req{216} into \req{211} and
considering the identity \req{214} we obtain the require.

In particular for $N=2$ we deduce the differential system
\begin{equation}\label{219}
\left\{%
  \begin{array}{ll}
   \dot{x} =\sum_{j=1}^S\tilde{\lambda}_j\{g_j,\,x\}+\tilde{\lambda}_{S+2}\{y,\,x\}\\
   \dot{y} =\sum_{j=1}^S\tilde{\lambda}_j\{g_j,\,y\}+\tilde{\lambda}_{S+2}\{x,\,y\}.
  \end{array}%
\right.
\end{equation}

 As usual we denote by
$\mathbb{R}[x]$  the ring of all real polynomials in the variables
$x\equiv(x_1,\,x_2,,,,,\,x_N)$. We consider the polynomial vector
field in $\mathbb{R}^N,$  with degree $n,$ i.e.,
\[
\begin{array}{cc}
\bold{v}=(v^1(x),...,v^N(x)),\quad{v^j(x)}\in{\mathbb{R}[x]},\quad j=1,2,..,N\\
n=max(deg(v^1(x)),...,deg(v^N(x)))
\end{array}%
\]

\textbf{ Definition 2.2}\,

{\it We say  that $\{g=0\}\subset\mathbb{R}^N$ is an invariant
algebraic hypersurface of the polynomial vector field $\textbf{v}$
of degree $n$ if there exists a polynomial $K\in\mathbb{R}[x]$
such that}
$$dg(\textbf{v})=K(x)\,g.$$
The polynomial $K$ at the degree at most $n-1$ is called the
cofactor of $g(x)=0.$

 \textbf{Definition 2.3}\,

 {\it A nonconstant (multivalued) function is said to be Darboux if
it is of the form}
$$f=\ln(\prod_{j=1}^Sg^{\sigma_j}_j(x)),\quad ,$$ where
$\sigma_j\in\mathbb{C},\,j=1,2,..,S$ are certain constants.

 \textbf{ Definition 2.4}\, {\it We shall say that the vector field
$\dot{\textbf{x}}=\textbf{v}(x)$ with
 invariant relations
\begin{equation}\label{220}
g_j(x)=0,\quad j=1,...,S>N,
\end{equation}
is integrable  if it admits $N-1$
independent  first integrals $f_1,\,f_2,\,..,f_{N-1},$ and
integrable in the Darboux sense if
$$g_1,\,g_2,.....,g_{S}$$
are polynomial functions and $f_1,\,f_2,...,f_{N-1}$ are Darboux
functions.}

 Darboux proved the following theorem.

\textbf{Darboux's theorem}

 {\it Si l'on connait $\frac{m(m+1)(m+2)....(m+n-1)}{n!}=M_n$
 int\'egreles particuli\'eres alg\'ebriques de syst\'eme

 $$\frac{dx_1}{L_1}=\frac{dx_2}{L_2}=....=\frac{dx_n}{L_n},\quad L_1,\,L_2,\,...,L_n\in\textsc{R}[x]$$
on pourra trouver lemultiplicateur du syst\'eme.

Si l'on connait $M_n+r$ int\'egrales particuli\'eres alg\'ebriques
du m\'eme syst\'eme, on pourra en determiner  le multiplicateur et
$r$ int\'egrales g\'en\'erales.

Si l'on connait $M_n+n-1=q$ int\'egrales particuli\'eres
alg\'ebriques $u_1,\,u_2,\,...,u_q$ on pourra effectuer
l'int\'egration compl\'ete. Les int\'egrales se pr\'esenteront
sous la forma suivante:}
\[
\begin{array}{cc}
u^{\alpha_1}_1u^{\alpha_2}_2...u^{\alpha_q}_q&=C_1\\
u^{\beta_1}_1u^{\beta_2}_2...u^{\beta_q}_q&=C_2\\
&\vdots\\
u^{\lambda_1}_1u^{\lambda_2}_2...u^{\lambda_q}_q&=C_{n-1}.
\end{array}%
\]

In \cite{Jou} the following result is proved.

\textbf{ Jounolou´s Theorem}

{\it Let $\textbf{v}$ be a polynomial vector field defined in
$\textsc{C}^N$ of degree $n>0.$ Then $\textbf{v}$ admits
$\frac{(n+N-1)!}{(n-1)!}+n$ irreducible invariant algebraic
hipersurface if and only if $\textbf{v}$ has a rational first
integral.}

\begin{prop}\, {\it The  vector field \req{211} with invariant relations
\req{220} is integrable if and only if the vector field the
Eruguin functions  are such that}

\begin{equation}\label{221}
 \Phi_j=\tilde{\lambda}\{f_1,\,f_2,\,..,f_{N-1},g_j\},\quad
k=1,2,..,S
\end{equation}
\end{prop}
where $\tilde{\lambda}$ is an arbitrary function.

\textbf{Proof.}\, Let us suppose that the vector field
$\textbf{v}$ is integrable, then it admits the representation
 \cite{Sad1},\,\cite{Ram3}
\[
\textbf{v}=\tilde\lambda \left|
\begin{array}{lllll}
df_1(\partial_{1})&df_1(\partial_{2})&\ldots&\vdots &df_1(\partial_{N})\\
df_2(\partial_{1})&df_2(\partial_{2})&\ldots& \vdots&df_2(\partial_{N})\\
\vdots         &\vdots  &\ldots                    &\vdots&\vdots\\
df_{N-1}(\partial_{1})&df_{N-1}(\partial_{2})&\ldots&\vdots &df_{N-1}(\partial_{N})\\
\partial_1&\partial_2&\ldots&\vdots&\partial_N
\end{array}%
\right| \equiv\tilde\lambda\,\{f_1,..., f_{N-1},*\},
\]
 where
$\tilde\lambda$ is an arbitrary function. Hence we obtain that
$$dg_j(\textbf{v})=\tilde\lambda\,\{f_1,..., f_{N-1},g_j\}$$
on the other hand from (2.2) we obtain that
$$dg_j(\textbf{v})=\Phi_j .$$ By compare both we deduce \req{221}.

 We obtain the reciprocity result as follows.

Let us suppose that \req{221} holds. Clearly that the condition
\req{214} holds identically in this case.

By inserting \req{221} in \req{211} we obtain

\[
\begin{array}{cc}
\textbf{v}(*)=\lambda(&\{f_1,f_2,f_3,...,f_{N-
1},g_1\}\{*,g_2,g_3,...,g_N\}+\\
...+&\{f_1,f_2,f_3,...,f_{N-1},g_N\}\{g_1,g_2,g_3,...,*\}).
\end{array}%
\]

In view of the identity \req{214} we deduce
$$\textbf{v}(*)=\{g_1,g_2,g_3,...,g_N\}\lambda\{f_1,f_2,f_3,...,f_{N-
1},*\}\equiv{\tilde\lambda\,\{f_1,..., f_{N-1},\,*\}}.$$ as a
consequence the vector field is  integrable.

\begin{cor} {\it Let us suppose that the system \req{2}8 is polynomial
of degree $n.$

Then it is Darboux integrable if and only if}
\begin{equation}\label{222}
\begin{array}{cc}
\tilde{\lambda}^\alpha_j=\kappa (x)\displaystyle\frac{\sigma_j\,f}{g_j}\quad j=1,2,..,S,\\
\tilde{\lambda}_{S+1}=\nu\{f,y\},...,
\,\tilde{\lambda}^\alpha_{S+2}=\nu\{x,f\}
\end{array}%
\end{equation}
where $\kappa ,\,\nu$ are arbitrary rational functions and $f$ is
a Darboux's function.
\end{cor}

\textbf{Proof.}

The proof follows from the fact that in view of
\req{214},\,\req{222} we obtain
\[
\left\{%
\begin{array}{cc}
dg_k(\textbf{v})=(\kappa + \nu)\, \{f, \,g_k\},\\
\\
\quad (\kappa + \nu) \{f, \,g_k\}|_{g_k(x)=0}=0,\quad k=1,2,..,S.
\end{array}%
\right.
\]
 {\textbf{ Example 2.2.}}

  We shall suppose that the given invariant relations of the differential system \req{211} are
the hyperplane $$x^j=0,\quad j=1,2,..,N.$$ We choose the Eruguin
functions as follows
$$\Phi_j=\Psi_j(x^j)\frac{\{\varphi_1,..., \varphi_{N-1},\,x^j\}}{\{\varphi_1,\varphi_2...,
\varphi_{N}\}},$$ where $\varphi_{kj}(x^k),\quad \Psi_j(0)=0,\quad
k,j=1,2,..,N$ and $\varphi_j,\,j=1,2,..,N$ are arbitrary
functions. Hence we obtain that this system takes the form

\begin{equation}\label{223}
\dot{x}^j=\Psi_j(x^j)\frac{\{\varphi_1,...,
\varphi_{N-1},\,x^j\}}{\{\varphi_1,\varphi_2...,
\varphi_{N}\}},\quad j=1,2,..,N
\end{equation}

We shall study the case when
\begin{equation}\label{224}
\left\{%
\begin{array}{cc}
 \{\varphi_1,\varphi_2..., \varphi_{N}\}\ne{0}\\
d\varphi_j=\sum_{k=1}^N\varphi_{kj}(x^k)dx^k, \quad j=1,2,..,N
\end{array}%
\right.
\end{equation}

The differential system \req{223},\,\req{224} is integrable.

 In fact, by considering that

\[\sum_{k=1}^N\frac{\varphi_{kj}(x^k)dx^k}{\Psi_{k}(x^k)}=\left\{
  \begin{array}{ll}
    dt&  \hbox{if $j=N$}, \\
    0 &  \hbox{if $j\ne{N}$}.
  \end{array}
\right.
\]

 we deduce the existence
of $N-1$ independents first integrals

\begin{equation}\label{225}
\left\{%
\begin{array}{cc}
f_{j}(x)\equiv\sum_{k=1}^N \int{\frac{\varphi_{kj}(x^k)}{\Psi
_k(x^k)}dx^k}=c_j,\quad j=1,2,...,N-1.
\end{array}%
\right.
\end{equation}

It is easy to show that the vector field $\textbf{v}(*)$ in this
case is such that
$$\textbf{v}(*)=g\,\displaystyle\frac{\{f_1,..., f_{N-1},\,*\}}{\{\varphi_1,
\varphi_2..., \varphi_{N}\}},\quad g=\prod_{k=1}^N\Psi_k(x^k).$$

 For the subcase when the invariant hyperplane are such that
 $$x^k-a_{m+k}=0,\quad k=1,2,..,N,\,m=1,2,...,M$$
 and

\begin{equation}\label{226}
\varphi_{kj}(x^k)=(x^k)^{j-1},\quad k,j=1,2,..,N\\
\Psi_k(x^k)=\prod_{m=1}^{M}(x^k-a_{m+k}),\end{equation}
 then the
first integrals \req{225} in this case take the form
\[
\left\{%
\begin{array}{cc}
f_j=\ln\prod_{k=1}^N{\prod_{m=1}^M(x^k-a_{m+k})^{\sigma^j_{k+m}}},\quad j=1,..,N\\
\\
\sigma^j_{m+k}=\displaystyle\frac{(a_{k+m})^j}{\prod_{l=1,\,l\ne{m}}^M(a_{k+l}-a_{k+m})}
\end{array}%
\right.
\]
as a consequence the system \req{221},\,\req{222},\,\req{224} is
Darboux integrable.

An interesting particular case is the following
\begin{equation}\label{227}
\dot{x}^j=x^j\prod_{m=1}^M(\frac{x^2_j}{m^2}-1)\frac{\{\varphi_1,...,
\varphi_{N-1},\,x^j\}}{\{\varphi_1,\varphi_2...,
\varphi_{N}\}}.\end{equation}
 Hence, by making $M\to{+\infty}$ we
deduce

\begin{equation}\label{228}
\dot{x}^j=x^j\prod_{m=1}^{+\infty}(\displaystyle\frac{x^2_j}{m^2}-1)\displaystyle\frac{\{\varphi_2,...,
\varphi_{N},\,x^j\}}{\{\varphi_1,\varphi_2...,
\varphi_{N}\}}\equiv{\sin{\pi
x^j}\displaystyle\frac{\{\varphi_2,...,
\varphi_{N},\,x^j\}}{\{\varphi_1,\varphi_2...,
\varphi_{N}\}}}.\end{equation}

We observe that the differential systems of the type \req{223}
appear in the theory of the St$\ddot{a}$ckel mechanical system
\cite{Sad1},\,\cite{Ram3}. With respect to this system we state
the following problem:

\textbf{Problem}

Determine the real constants $K_1,\,K_2,\,...,K_{N-1},\,L$ in such
a way that the  hyperplane

\begin{equation}\label{229}
x^N=\sum_{j=1}^{N-1}K_jx^j+L
\end{equation}

 is invariant of the  system  \req{221},\,\req{224}.

We solve this problem for the case when
\begin{equation}\label{230}
\left\{%
\begin{array}{cc}
 N=2,\quad
\varphi_2=\frac{a}{2}(x^2+y^2)+xy,\\
\Psi_1(x)=-\lambda\prod_{j=1}^M(\displaystyle\frac{x^2}{a_j^2}-1),\quad
\Psi_2(y)=\lambda\prod_{j=1}^M(\displaystyle\frac{y^2}{a_j^2}-1)
\end{array}%
\right.
\end{equation}
 where $a,a_1,...,a_M$ are real constants and $\lambda$ is an arbitrary
function.

 The system \req{221} in this case  takes the form
\begin{equation}\label{231}
\left\{%
\begin{array}{cc}
\dot{x}=\lambda\,(x+ay)\prod_{j=1}^M(\displaystyle\frac{x^2}{a_j^2}-1)\\
\dot{y}=\lambda\,(ax+y)\prod_{j=1}^M(\displaystyle\frac{y^2}{a_j^2}-1)
\end{array}%
\right.
\end{equation}

 We require to determine the real values of the constants $K$ and $L$ in such a way
 \begin{equation}\label{230}y=Kx+L\end{equation}
 is an invariant straight line of \req{231}.

 Clearly that the parameter $K$ must be satisfies the relation
 $$K^{2M}+aK^{2M-1}-aK-1=0.$$
Hence we obtain that  $K_1=1,\,K_2=-1$ satisfies this relation.
For $M>2$ there exist at most four real values of $K$ which
satisfy this equation.

 By using the
algebraic computer packages we can solve the stated problem. In
particular for the cubic and quintic  system in which

$$\lambda=(\prod_{j=1}^Ma^2_j)^{-1}$$
and
\[
\begin{array}{cc}
a=\sqrt{5},\quad a_1=1,\quad M=1\\
a=\sqrt{5},\quad a_1=1,\quad a_2=\sqrt{5}-2,\quad M=2
\end{array}%
\]

It is easy to show that in this case we obtain respectively

$$K_1=1,\quad L_1=0\quad  K_2=-1,\quad L_2=0$$
and
\[
\begin{array}{cccc}
 K_1 = 1,\quad L_1 = 0, \quad  K_2 = -1,\quad L_2 = 0,\\
 K_3=-\frac{1}{2}-\frac{\sqrt{5}}{2},\quad L_3 = -\frac{1}{2}+\frac{\sqrt{5}}{2}\\
K_4 = \frac{1}{2}-\frac{\sqrt{5}}{2},\quad L_4 =
-\frac{3}{2}+\frac{\sqrt{5}}{2},\\
K_3 = -1/2-\frac{\sqrt{5}}{2},\quad L_5 = 1/2-\frac{\sqrt{5}}{2},\\
K_4 = 1/2-\frac{\sqrt{5}}{2},\quad L_6 =
\frac{3}{2}-\frac{\sqrt{5}}{2}
\end{array}%
\]

 The  quintic polynomial system in
this case was constructed in \cite{Art}

\[
\left\{%
\begin{array}{cc}
\dot{x}=(x+\sqrt{5}y)(x^2-1)(x^2-(\sqrt{5}-2)^2)\\
\dot{y}=(\sqrt{5}x+y)(y^2-1)(y^2-(\sqrt{5}-2)^2)
\end{array}%
\right.
\]
 and admits 14
straight lines.

The Eruguin functions in this case are:
\[
\left\{%
\begin{array}{cc}
\Phi_m=(x+\sqrt{5}\,y)(x^2-1)(x^2-(\sqrt{5}-2)^2),\\
\Phi_{4+m}=(\sqrt{5}\,x+y)(y^2-1)(y^2-(\sqrt{5}-2)^2),\quad
m=1,2,3,4
\end{array}%
\right.
\]

 hence  the constructed system  is not
Darboux integrable.

  For the case when in \req{231} $$a_j=j,\quad a=0$$ we obtain the
system

\begin{equation}\label{229}
\left\{%
\begin{array}{cc}
\dot{x}=x\prod_{j=1}^M(\displaystyle\frac{x^2}{j^2}-1)\\
\dot{y}=y\prod_{j=1}^M(\displaystyle\frac{y^2}{j^2}-1)
\end{array}%
\right.
\end{equation}
 In this case the system
admits the following invariant straight line
\[
\left\{%
\begin{array}{cc}
x=j,\quad x=-j,\quad j=1,2,..,M\\
y=j,\quad y=-j\\
y=x,\quad y=-x
\end{array}%
\right.
\]
 Clearly, the system \req{229}is Darboux integrable

 If $M=+\infty$ then the system \req{229}
takes the form

\[
\left\{%
\begin{array}{cc}
\dot{x}=x\prod_{j=1}^\infty(\frac{x^2}{j^2}-1)=\sin{\pi x}\\
\dot{y}=y\prod_{j=1}^\infty(\frac{y^2}{j^2}-1)=\sin{\pi
y}\end{array}%
\right.
\]

 for which the infinity numbers of the straight
lines

\[
\left\{%
\begin{array}{cc}
x=j,\,x=-j,\quad j=1,2,..,+\infty\\
y=j,\quad y=-j\quad j=1,2,..,+\infty\\
y=x+2m,\quad y=-x+2m,\quad m\in\textsc{Z}
\end{array}%
\right.
\]
are its invariant.

The problem of the determination of the upper bound for the
maximum number of the invariant straight lines ($L(n)$)  for the
polynomial system is an open problem.

It is easy to show that \cite{Ram2}

\[L(n)\geq
\left\{
  \begin{array}{ll}
    {{2n+1}} &  \hbox{if $n$ is even}, \\
   {2n+2}&  \hbox{if $n$ is odd}.
  \end{array}%
  \right.
\]

 There exist the following conjecture

 {\textbf{Conjecture }\cite{Art}

$$L(n)\leq{3n-1}$$

This upper bound is reached in particular for $n=2,3,4,5.$

\section{ Inverse approach for the planar vector fields}

 In this section we analyze the
Eruguin-Galliulin theory developed in the above section for the
case when $N=2.$

The differential system \req{211} in this case take the form
\begin{equation}\label{31}
\left\{%
\begin{array}{cc}
\dot{x}&=\Phi_1\{x,g_2\}+\Phi_2\{g_1,x\}=P(x,y)\\
\dot{y}&=\Phi_1\{y,g_2\}+\Phi_2\{g_1,y\}=Q(x,y),
\end{array}%
\right.
\end{equation}

 we set
$\textbf{v}=(P,Q).$

For the case when this equations admit the subsidiary invariant
curves
$$g_j(x,y)=0,\quad j=3,4,...,S,$$
The Eruguin functions must be satisfy the relations

\begin{equation}\label{32}
\Phi_1\{g_j,g_2\}+\Phi_2\{g_1,g_j\}+\Phi_j\{g_2,g_1\}=0\quad
j=1,2,...S.
\end{equation}

Hence the  Eruguin functions  $\Phi_m$ can be determine as follows
\begin{equation}\label{33}
\left\{%
\begin{array}{cc}
\Phi_m=\sum_{j=1}^S\tilde{\lambda}_j(x)\{g_j,\,g_m\}\prod_{k=1,\,k\ne{j}}^Sg_k
+(\tilde{\lambda}_{S+1}\{x,\,g_m\}+\tilde{\lambda}_{S+2}\{\,g_m,\,y\})g,\quad
\\
g=\prod_{j=1}^Sg_j
\end{array}%
\right.
\end{equation}
 where
$\tilde{\lambda}_1,\,\tilde{\lambda}_2,...,\tilde{\lambda}_{S+2},$
are arbitrary functions.

The veracity of this representation we obtain  by inserting
\req{33} into \req{32} and by considering the identity \req{214}
which in this case takes the form
\begin{equation}\label{34}
 \{g_k,\,g_j\}\{g_i,\,g_m\}+ \{g_k,\,g_m\}\{g_j,\,g_i\}+
\{g_i,\,g_k\}\{g_j,\,g_m\}\equiv{0},
\end{equation}

\textbf{ Corollary 3.1}\,

{\it The differential equations \req{31},\,\, \req{33} can be
rewritten as follows}

\begin{equation}\label{35}
\left\{%
\begin{array}{cc}
 \dot{x}=g(x,y)(
\sum_{j=1}^S\tilde{\lambda}_j(x,y)\frac{\{g_j,x\}}{g_j}+\tilde{\lambda}_{S+2})=P(x,y)\\
\dot{y}=g(x,y)(
\sum_{j=1}^S\tilde{\lambda}_j(x,y)\frac{\{g_j,y\}}{g_j}-\tilde{\lambda}_{S+1})=Q(x,y)
\end{array}%
\right.
\end{equation}

 We deduce the proof by inserting \req{33} into \req{31} by using the identity
\req{34}.

From these relations we obtain the following proposition

 \textbf{Proposition 3.1}\,

 {\it Let $g_j(x,y)=0,\,j=1,2,...,S$ are the irreducible algebraic
curves, then the polynomial differential system \req{35} admits
the Darboux first integral

$$f(x,y)=\ln(\prod_{j=1}^Sg^{\sigma_j}_j)$$
if and only if in \req{33}}
 $$\tilde{\lambda}_j=\nu_0 \sigma_j=constants,\quad \tilde{\lambda}_{S+1}=
 \nu \{f,y\},\quad \tilde{\lambda}_{S+2}=\nu \{x,f,\},$$
where $\nu_0,\nu$ are arbitrary rational functions and
$\sigma_j=constants,\,j=1,2,...,S.$

The proof is easy to obtain from corollary 2.2.

The system \req{35} in this case takes the form

\begin{equation}\label{35}
\left\{%
\begin{array}{cc}
\dot{x}=V\{x,\,F\}\\
\dot{y}=V\{y,\,F\}
\end{array}%
\right.
\end{equation}
 where
$$F=\prod_{j=1}^Sg^{\sigma_j}_j,\quad V=g(\frac{\nu_0}{F}+\nu ).$$

 We illustrate the above results in the following
concrete cases.

 \textbf{Example 3.1}

In this section we give the results exposed in \cite{Sad1},
related with the construction the planar polynomial vector field
with invariant circumferences:
$$g_j(x,y)\equiv{(x-a_j)^2+(y-b_j)^2-r^2_j}=0,\quad j=1,2,...,S$$

The system \req{35} under the restrictions
$$\tilde{\lambda}_{S+1}=\tilde{\lambda}_{S+2}=0$$
takes the form
\begin{equation}\label{36}
\left\{%
\begin{array}{cc}
\dot{x}=-\sum_{j=1}^S\tilde{\lambda}_j(y-b_j)\prod_{m=1,\,m\ne{j}}^Sg_m\equiv{P(x,y)}\\
\dot{y}=\sum_{j=1}^S\tilde{\lambda}_j(x-a_j)\prod_{m=1,\,m\ne{j}}^Sg_m\equiv{Q(x,y)},
\end{array}%
\right.
\end{equation}

or, what is the same,
\[
\left\{%
\begin{array}{cc}
\dot{x}=-y(x^2+y^2)^{S-1}\sum_{j=1}^S\tilde{\lambda}_j-(x^2+y^2)^{S-1}\sum_{j=1}^S\tilde{\lambda}_jb_j+..\equiv{P(x,y)}\\
\dot{y}=x(x^2+y^2)^{S-1}\sum_{j=1}^S\tilde{\lambda}_j+(x^2+y^2)^{S-1}\sum_{j=1}^S\tilde{\lambda}_ja_j+...\equiv{Q(x,y)}.
\end{array}%
\right.
\]

Now we determine the arbitrary functions
$\lambda_1,\,\lambda_2,\,...,\lambda_S$ in such a way that the
above vector field is polynomial of the fixed degree $n.$

\textbf{Corollary 3.2}\,

{\it Let  us suppose that
\[
\left\{%
\begin{array}{cc}
\tilde{\lambda}_j\in{\Bbb{R}[x,y]},\quad j=1,2,..,S\\
S_1=deg\sum_{j=1}^S\tilde{\lambda}_j,\\
S_2=max(deg(\sum_{j=1}^S\tilde{\lambda}_jb_j),deg(\sum_{j=1}^S\tilde{\lambda}_ja_j)
\end{array}%
\right.
\]

then
 $$n=max(deg(P),deg(Q))\leq 2S-1+S_1,$$ if\quad
$$\sum_{j=1}^S\tilde{\lambda}_j\ne{0},$$ and
$$n=max(deg(P),deg(Q))\leq 2S-2+S_2,$$ if
$$\sum_{j=1}^S\tilde{\lambda}_j={0}.$$}

From this results we obtain the proof of the following result

\textbf{ Proposition 3.2}\,

{\it Every configuration of the circumferences in the plane is
realizable by a polynomial of the degree at most $2S+S_1-1$ or
$2S+S_2-2$ where $S_j,\,j=1,2$  are the degree of the polynomials
introduced above.}

 In a paper \cite{Llib4} the authors proved that
every configuration of cycles on the plane is realizable (up to
homeomorphism) by a polynomial of the degree at most $2(m+r)-1,$
where $m$ is the number of cycles and $r$ is the number of primary
cycles  (a cycle $C$ is primary if there are no other cycles
contained in the bounded region limited by $C$).

It is interesting to observe that the upper bound for the degree
of the constructed vector field is  independent from whether its
cycles are primary or not.

Now we shall study the case when the  circumferences form two
nests with the centers at  the points $(0,0)$ and $(a,0)$
respectively, hence

\[
\left\{%
\begin{array}{cc}
g_j(x,y)\equiv{x^2+y^2-r^2_j}=0,\quad j=1,2,...,l_1,\\
g_j(x,y)\equiv{(x-a)^2+y^2-r^2_j}=0,\quad j=l_1,...,S
\end{array}%
\right.
\]

Clearly,  for this case the Eruguin functions are such that
\[
\left\{%
\begin{array}{cc}
\Phi_1=\Phi_3=.....=\Phi_{l_1}=4ay\sum_{j=l_1+1}^{S}\tilde{\lambda}_j\prod_{m=1,\,m\ne{2}}^{l_1}g_m\\
\Phi_2\equiv{\Phi_{l_1+1}}=\Phi_{l_1+2}=.....=\Phi_{S}=-4ay\sum_{j=1}^{l_1}\tilde{\lambda}_jg_2\prod_{m=l_1+1}^Sg_m.
\end{array}%
\right.
\]

 The
differential system \req{36} takes the form
\[
\left\{%
\begin{array}{cc}
\dot{x}=-2(\sum_{j=1 }^{l_1}\tilde{\lambda}_j\prod_{k=l_1+1}^Sg_k+
\sum_{j=l_1+1}^S\tilde{\lambda}_j\prod_{k=l}^{l_1}g_k)y,\\
\dot{y}=2(\sum_{j=1 }^{l_1}\tilde{\lambda}_j\prod_{k=l_1+1}^Sg_k+
\sum_{j=l_1+1}^S\tilde{\lambda}_j\prod_{k=l}^{l_1}g_k)x-2a\sum_{j=l_1+1}^S\tilde{\lambda}_j\prod_{k=l}^{l_1}g_k
.\end{array}%
\right.
\]

In particular for the case when
\[
\left\{%
\begin{array}{cc}
n=2l+1=S+1,\quad l_1=l,\\
-2\sum_{j=1 }^{l}\tilde{\lambda}_j=x-a+y,\\
-2\sum_{j=l+1 }^{S}\tilde{\lambda}_j=x+y,
\end{array}%
\right.
\]

 we obtain the
vector field constructed in \cite{Sad1}.

By designating by
 $F_{a}(x,y),\,F_0(x,y)$ the following polynomials
\[
\left\{%
\begin{array}{cc}
F_{a}(x,y)&=(x+y-a)\prod_{j=2}^{l+1}((x-a)^2+y^2-r^2_j)\quad{l\geq{1}},\\
F_0(x,y)&=F_{a}(x,y)|_{a=0}.
\end{array}%
\right.
\]
we can deduce that the above vector field takes the form:
\[
\left\{%
\begin{array}{cc}
\dot{x}&=\Big(F_0(x,y)-F_{a}(x,y)\Big)\,y=P(x,y)\\
\dot{y}&=-\Big(F_0(x,y)-F_{a}(x,y)\Big)\,x+aF_0(x,y)=Q(x,y).\end{array}%
\right.
\]

This system has the following properties:

 \quad 1) has
only 3 critical points in the finite plane
 $\textsc{R}^2$
$$(0,0),\quad{(\displaystyle\frac{a}{2},\,0),} \quad{({a},\,0)}. $$

\quad 2) the Liapunov quantities $\sigma $ and $\Delta$ for the
system are :

i)\[
\left\{%
\begin{array}{cc}
\sigma{(0,0)}=\sigma{(a,0)}\\
\Delta{(0,0)}=\Delta{(a,0)}.\end{array}%
\right.
\]

ii))\[
\left\{%
\begin{array}{cc}
\sigma{(0,0)}=(-1)^la\prod_{j=1}^{l}r^2_j\\
\Delta{(0,0)}= a^2\prod_{j=1}^l
l(a^2-r^2_j)\Big(\prod_{j=1}^l(a^2-r^2_j)-(-
1)^l\prod_{j=1}^lr^2_j\Big)\\
\sigma^2{(0,0)}-4\Delta{(0,0)}=
a^2\Big(\big(2\prod_{j=1}^l(a^2-r^2_j)-(-1)^l\prod_{j=1}^lr^2_j\big)^2-
8\big(\prod_{j=1}^l(a^2-r_j^2)\big)^2\Big).
\end{array}%
\right.
\]

\[
\left\{%
\begin{array}{cc}
\sigma{(\displaystyle\frac{a}{2},0)}=a\prod_{j=1}^l\big((\displaystyle\frac{a}{2})^2-r^2_j\big)\\
\Delta{(\displaystyle\frac{a}{2},0)}=
-\displaystyle\frac{a^4}{2}\prod_{j=1}^{l}\big((\displaystyle\frac{a}{2})^2-r^2_j\big)
 \sum_{l=1}^l\prod_{j=1,\,j\ne{l}}^l\Big((\displaystyle\frac{a}{2})^2-r^2_j\Big).
\end{array}%
\right.
\]

The circumferences do not intersect if
 $$r_j<a/2,\quad j=1,...,l,$$
so
$$\Delta{(\displaystyle\frac{a}{2},0)}<0, $$
and, as a consequence the critical point
 $(\displaystyle\frac{a}{2},0)$ is a saddle.

It is evident that the other critical points are the stability or
non stability foci depending on whether $k$ is odd or even.

Hence we obtain that the constructed polynomial vector field of
degree $n=S+1$ admits $S=2l$ invariant circumferences.

The proposition 3.1 we illustrate in the next two examples.

\textbf{ Example 3.2}\,

The particular case of the Lienard equation
$$\ddot{x}-\displaystyle\frac{d}{dx}h(x)\dot{x}-\alpha h(x)\displaystyle\frac{d}{dx}h(x)=0$$

or, what is the same,
\[
\left\{%
\begin{array}{cc}
\dot{x}&=y+h(x)\\
\dot{y}&=\alpha h(x)\displaystyle\frac{d}{dx}h(x)\end{array}%
\right.
\]

 is
Darboux integrable.

In fact, the first integral $F$ in this case is the following
$$F(x,y)=g^{\sigma_1}_1g^{\sigma_2}_2$$
where $g_1,\,g_2,\,\sigma_1,\,\sigma_2$ are such that
\[
\left\{%
\begin{array}{cc}
g_1(x,y)=y+\displaystyle\frac{1+\sqrt{4\alpha+1}}{2}h(x)\\
g_2(x,y)=y+\displaystyle\frac{1-\sqrt{4\alpha+1}}{2}h(x)\\
\sigma_1=-\displaystyle\frac{1-\sqrt{4\alpha+1}}{2\sqrt{4\alpha+1}},\quad\sigma_2=\displaystyle\frac{1+
\sqrt{4\alpha+1}}{2\sqrt{4\alpha+1}}
\end{array}%
\right.
\]

It is easy to show that in this case the Eruguin functions are
$$\Phi_1=\displaystyle\frac{1+\sqrt{4\alpha+1}}{2\sqrt{4\alpha+1}}g_1,\quad \Phi_2=\displaystyle\frac{1-\sqrt{4\alpha+1}}{2\sqrt{4\alpha+1}}g_2$$

Clearly, if $4\alpha+1<\,0$ then the function $F$ takes the form
$$F(x,y)=(y^2+h(x)y-\alpha
 h^2(x))\exp\Big(\sqrt{-\alpha-\displaystyle\frac{1}{4}}arctan{\displaystyle\frac{\sqrt{-\alpha-
 \displaystyle\frac{1}{4}}h(x)}{y+\displaystyle\frac{1}{2}h(x)}}\Big)$$

It is interesting to observe that if the function $h$ admits the
following development
$$h(x)=x+a_2x^2+...$$
then the origin of the given system is a focus.

\textbf{Example 3.3}\,

The differential equation

\begin{equation}\label{37}
 \dot{z}=i\,(a_{10}z+a_{01}\bar{z}+\sum_{j+k=3}a_{jk}z^j\,\bar{z}^{k})
\end{equation}
is Darboux integrable,\, where $a=(a_{jk}),\, j,k=0,2,3$ are real
constants matrix and
$$  z=x+iy,\quad \bar{z}=x-i\,y$$ are the complex
coordinate in the plane $\mathbb{R}^2.$

In fact, the equations \req{37} are equivalent to the cubic planar
system
\[
\left\{%
\begin{array}{cc}
\dot{x}=y(a_{01}-a_{10}+(a_{12}-a_{21}+3(a_{03}-a_{30})\,x^2+
(a_{12}-a_{21}+a_{30}-a_{03})\,y^2)\\
\dot{y}= x(a_{01}+a_{10}+(a_{12}+a_{21}+a_{03}+a_{30})\,x^2+
(a_{12}+a_{21}-3(a_{30}+a_{03})y^2).
\end{array}%
\right.
\]

 Hence, by introducing the correspondent notations we obtain the system
\begin{equation}\label{38}
\left\{%
\begin{array}{cc}
\dot{x}=y(a+b\,x^2+c\,y^2)\\
\dot{y}=x(\alpha+\beta\,x^2+\gamma\,y^2).
\end{array}%
\right.
\end{equation}
 We shall analyze the
case when $c\ne{0}.$

Let $g_1,\,g_2$ are the functions:
$$g_j(x,y)=\nu_j(x^2-\lambda_0)-y^2+\lambda_1,\quad j=1,2$$
where $\lambda_0,\,\lambda_2,\,\nu_1,\,\nu_2$ are constants:
\[
\left\{%
\begin{array}{cc}
&\lambda_0=\displaystyle\frac{\gamma\,a-\alpha\,c}{{b\gamma-c\beta}},\quad \lambda_1=\displaystyle\frac{\alpha\,b-\beta\,a}{{b\gamma-c\beta}}, \\
&\nu_1=\displaystyle\frac{\gamma-b}{2c}+\sqrt{(\displaystyle\frac{\gamma-b}{2c})^2+\displaystyle\frac{\beta}{c}},\quad
\nu_2=\displaystyle\frac{\gamma-b}{2c}-\sqrt{(\displaystyle\frac{\gamma-b}{2c})^2+\displaystyle\frac{\beta}{c}}\\
&\nu_1-\nu_2\ne{0},
\end{array}%
\right.
\]
then the following relations hold

\begin{equation}\label{38}
\left\{%
\begin{array}{cc}
dg_j(\textbf{v})=2x\,y\,(\gamma-\nu_j\,c)g_j,\quad j=1,2\\
\{g_1,\,g_2\}=4x\,y\,(\nu_1-\nu_2).
\end{array}%
\right.
\end{equation}

The proof is easy to obtain after some calculations.

 The given vector field is Darboux integrable with $F:$
$$F(x,y)=\displaystyle\frac{(\nu_1(x^2-\lambda_0)-y^2+\lambda_1)^{b+\nu_1c}}{(\nu_2(x^2-\lambda_0)-y^2+\lambda_1)^{b+\nu_2c}},$$
here we use the relation
$$c(\nu_1+\nu_2)=\gamma-b.$$

Now we shall  study the case  when $\nu_1,\,\nu_2$ are complex
numbers.

By introducing the notations
$$\gamma-b=2cq,\quad \gamma+b=2cr$$ we obtain that

$$\nu_1=q+i\,p,\quad \nu_2=q-i\,p ,\quad p^2=-4\beta\,c-q^2,\,p\,>0$$
The system (3.9) takes then the form (we put $c=1$)
\begin{equation}\label{38}
\left\{%
\begin{array}{cc}
\dot{x}=y(a+(r-q) \,x^2+y^2)\\
\dot{y}=x(\alpha-(p^2+q^2) \,x^2+(r+q)y^2).
\end{array}%
\right.
\end{equation}
 By considering that
in this case
$$g_1(x,y)=q(x^2-\lambda_0)-y^2+\lambda_1+ip(x^2-\lambda_0),$$
 we obtain that the first integral $F$ takes the
form
$$F(x,y)=\Big((y^2-\lambda_1-q(x^2-\lambda_0))^2+p^2(x^2-\lambda_0)^2\Big)
exp\Big({2r\,arctg\displaystyle\frac{p(x^2-\lambda_0)}
{y^2-\lambda_1-q(x^2-\lambda_0)}}\Big).$$

\section{ Planar differential system with one invariant algebraic
curve}

In this section, by applying the results of the section 3, we
construct the analytic planar vector field
\begin{equation}\label{41}
\left\{%
\begin{array}{cc}
\dot{x}=P(x,y)\\
\dot{y}=Q(x,y),\end{array}%
\right.
\end{equation}
 where $P$ and $Q$
are analytic functions on the region\, $G\subset\mathbb{R}^2,$
 from a given set of trajectories:

\begin{equation}\label{42}
\left\{%
\begin{array}{cc}
g_j(x,y)=y-y_j(x)=0,\quad j=1,2,\ldots,S\geq{2}\\
\{g_1,\,g_2\}=y^{'}_2(x)-y^{'}_1(x)\ne{0},\\
\prod_{j=1}^S y_j(x)\ne{0}
\end{array}%
\right.
\end{equation}
where $\displaystyle\frac{d{y_j}}{dx}=y^{'}_j$ and
$y_1,\,y_2,\,\ldots,y_S\in{C^{r}}(G\subset{\textsc{R}}),$
$r\geq{1}.$

 We shall study the particular case when $y_1,\,y_2,\,...y_s$ are
 solutions of the equation
 $$g(x,y_j)=0,\quad j=1,2,..,S$$
where $g(x,y)=0$ is an algebraic irreducible curve.

 By considering
\req{31} we obtain that the require system \req{41},\,\req{42}:
\begin{equation}\label{43}
\left\{%
\begin{array}{cc}
\dot{x}&=\Phi_1(x,y)-\Phi_2(x,y)\\
\dot{y}&=\Phi_1(x,y)y^{'}_2(x)-\Phi_2(x,y)y^{'}_1(x)
\end{array}%
\right.
\end{equation}

The condition \req{32}, on the existence of the complementary
partial integrals, in this case take the form

\begin{equation}\label{44}
\Phi_1({y_2}(x)- {y_j}(x))){'}+
  \Phi_2({y_j}(x)- {y_1}(x)){'}+
\Phi_j({y_1}(x)- {y_2}(x)){'}=0, \quad j=3,4\ldots,S.
\end{equation}
The given set of differential equations \req{35} in this case can
be rewritten as follows

\begin{equation}\label{45}
\left\{%
\begin{array}{cc}
\dot{x}=\sum_{j=1}^S\tilde{\lambda}_j\prod_{m\ne{j}}(y-y_m)+\tilde{\lambda}_{S+2}g=P(x,y)\\
\dot{y}=\sum_{j=1}^S\tilde{\lambda}_jy^{'}_j\prod_{m\ne{j}}(y-y_m)-\tilde{\lambda}_{S+1}g=Q(x,y)
\end{array}%
\right.
\end{equation}

or, what is the same,
\begin{equation}\label{46}
\left\{%
\begin{array}{cc}
\dot{x}=y^S\tilde{\lambda}_{S+2}+y^{S-1}(\sum_{j=1}^S\tilde{\lambda}_j-\tilde{\lambda}_{S+2}\sum_{j=1}^Sy_j)+\\
y^{S-2}(\sum_{j=1}^S\tilde{\lambda}_jy_j-
\sum_{j=1}^Sy_j\sum_{k=1}^S\tilde{\lambda}_k)+..)+...\\
\dot{y}=-y^S\tilde{\lambda}_{S+1}+y^{S-1}(\sum_{j=1}^S\tilde{\lambda}_jy^{'}_j-\tilde{\lambda}_{S+1}\sum_{j=1}^Sy_j)+\\
y^{S-2}(\sum_{j=1}^S\tilde{\lambda}_jy^{'}_j-
\sum_{j=1}^Sy_j\sum_{k=1}^S\tilde{\lambda}_ky_ky^{'}_k...)+....
\end{array}%
\right.
\end{equation}

we set $\textbf{v}=(P,Q).$

\begin{prop}

{\it Let us suppose that the arbitrary functions
$\tilde{\lambda}_1,\,...\tilde{\lambda}_S,\,\tilde{\lambda}_{S+1},\,\tilde{\lambda}_{S+2}$
are such that
\begin{equation}\label{47}
\left\{%
\begin{array}{cc}
\tilde{\lambda}_{S+1}=-q_0(x),\quad\tilde{\lambda}_{S+2}=p_0(x)\\
\tilde{\lambda}_j=\displaystyle\frac{\triangle_j}{\triangle_0}(\sum_{k=0}^Sp_k(x)y^{S-k}_j(x))\\
\tilde{\lambda}_jy^{'}_j=\displaystyle\frac{\triangle_j}{\triangle_0}(\sum_{k=0}^Sq_k(x)y^{S-k}_j(x)),\end{array}%
\right.
\end{equation}
where $p_k,\,q_k$ are continuous functions on
$D\subset\mathbb{R}.$
\[\triangle_0=\left|
\begin{array}{llllll}
1&1&\ldots&1&\vdots&1\\
y_1&y_2&\ldots&y_j&\vdots&y_S\\
\vdots&\vdots&\ldots&\vdots&\vdots\\
y^{S-1}_1&y^{S-1}_2&\ldots&y^{S-1}_j&\vdots&y^{S-1}_S
\end{array}
\right|
\]

\[\triangle_j=\left|
\begin{array}{lllllll}
1&1&\ldots&1&1&\vdots&1\\
y_1&y_2&\ldots&y_{j-1}&y_{j+1}&\vdots&y_S\\
\vdots&\vdots&\ldots&\vdots&\vdots\\
y^{S-2}_1&y^{S-2}_2&\ldots&y^{S-2}_{j-1}&y^{S-2}_{j+1}&\vdots
&y^{S-2}_S ,
\end{array}
\right|
\]
 Then the differential system \req{46} takes the form}

\begin{equation}\label{48}
\left\{%
\begin{array}{cc}
\dot{x}=&p_0(x)y^{S}+....+p_{S}(x)=P(x,y)\\
\dot{y}=& q_0(x)y^{S} +.... +q_{S}(x)=Q(x,y)
\end{array}%
\right.
\end{equation}
\end{prop}
and we set $\textbf{v}=(P,Q).$

 We shall
study the case when $p_0,\,p_1,\,...,p_S,\,q_0,\,q_1,\,...q_S$ are
polynomials on the variable $x$
 and such that $\textbf{v}$
represented a polynomial vector field of degree
$$n=max(degP,\,degQ)$$

\begin{cor}\, {\it Let g be a irreducible polynomial on the variables
$x$ and $y$:
\begin{equation}\label{49}
g=a_0(x)\prod_{j=1}^S(y-y_j(x))=\sum_{j=0}^Sa_j(x)y^{S-j},
\end{equation}

where\begin{equation}\label{410}
\left\{%
\begin{array}{cc}
a_1=-a_0(x)\sum_{j=1}^Sy_j(x)\\
a_2=a_0(x)\prod_{j<k}y_j(x)y_k(x)\\
&\vdots\\
a_S=(-1)^Sa_0(x)\prod_{j=1}^Sy_j(x)
\end{array}%
\right.
\end{equation}
then
$$dg(\textbf{v})=(\sum_{j=1}^S\displaystyle\frac{\Phi_j}{g_j})g=K(x)g$$}
\end{cor}
where $K$ is the cofactor of $g=0.$

\begin{prop}

{\it Let the curve \req{49}\req{410}) is invariant of the non zero
polynomial system of degree $n:$
\begin{equation}\label{411}
\left\{%
\begin{array}{cc}
\dot{x}=r_1(x)y^{n-1}+....+r_{n}(x)=P(x,y)\\
\dot{y}= q_0(x)y^{n} +.... +q_{n}(x)=Q(x,y)
\end{array}%
\right.
\end{equation}

where $r_j(x),\,q_j(x),\,j=0,1,..,n$ are polynomials of degree $j$
in the variable $x,$ then
$$S\leq{2n}.$$}
\end{prop}
 Proof, al absurd, let us suppose that $S=2n+1,$ then
from (4.9),\,(4.7),\, (4.11) we obtain that
\[
\left\{%
\begin{array}{cc}
p_{n+1}(x)=0,\\
p_j(x)=0,\\
q_j(x)=0,\quad j=0,1,..n,
\end{array}%
\right.
\]

 on the other hands from \req{47},\,\req{48} we
obtain that $$\tilde{\lambda}_j(x)=0,\quad j=1,2,...,2n+2$$ hence
the vector field is a zero vector field. Contradiction.

\section{Quadratic system  with one invariant algebraic curve}

In this section we shall study the case in which the vector field
\req{45} is quadratic i.e.,

\[
\left\{%
\begin{array}{cc}
 \dot{x}=p_{S-2}y^{2}+p_{S-1}(x)y+p_S(x)=P(x,y)\\
\dot{y}= q_{S-2}y^{2}+q_{S-1}(x)y+q_S(x)=Q(x,y)
\end{array}%
\right.
\]

 where $\max{(degP,\,degQ)}=2$ and  $q_{S-j},\,p_{S-j},\,j=0,1,2$ are polynomials in the
 variable $x.$ Below, for simplicity we shall denote this system as
 follows

 \begin{equation}\label{51}
\left\{%
\begin{array}{cc}
 \dot{x}=&p_{0}y^{2}+p_{1}(x)y+p_2(x)=P(x,y)\\
\dot{y}=& q_{0}y^{2}+q_{1}(x)y+q_2(x)=Q(x,y)
\end{array}%
\right.
\end{equation}

 First we prove the following general results related with the system
 \req{51}.

\begin{prop}
 {\it Let us suppose that \req{51} is such
that}
\begin{equation}\label{52}
dg(\textbf{v})=(\alpha_0y+\alpha x+\beta )g,
\end{equation}
\end{prop}
 where $g$ is given in
the formula \req{49},\req{410} and  $ \alpha_0,\,\alpha ,\,\beta $
are real constants and $p_j,\,q_j$ are polynomials of degree $j$
in the variable $x:$
\begin{equation}\label{53}
\left\{%
\begin{array}{cc}
p_j=\sum_{k=0}^jp_{jk}x^k,\\
q_j=\sum_{k=0}^jq_{jk}x^k,\quad j=0,1,2.
\end{array}%
\right.
\end{equation}
Then,

  \qquad If $p_0\ne{0}$ hence
\begin{equation}\label{54}
\left\{%
\begin{array}{cc}
\max (deg \,a_j(x))\leq{j},\\
\max (deg \,g)\leq{S}.
\end{array}%
\right.
\end{equation}

\qquad If
\begin{equation}\label{55}
\left\{%
\begin{array}{cc}
p_0={0},\quad p_{11}\ne{0},\quad \\
\alpha_0=(Sk+m)p_{11},\quad q_0=kp_{11} \\
  Sk+m\in{\textsc{N}}
\end{array}%
\right.
\end{equation}

hence}
\begin{equation}\label{56}
\left\{%
\begin{array}{cc}
\max (deg\, a_j(x))\leq{kj+m},\quad j=1,2,..,S\\
\max (deg\, g)\leq{Sk+m}.)
\end{array}%
\right.
\end{equation}

 In fact, from \req{52} we
deduced the differential system
\begin{equation}\label{57}
\left\{%
\begin{array}{cc}
A\cdot\displaystyle\frac{d\bold{a}}{dx}=B\cdot\bold{a}\\
p_0\displaystyle\frac{d{a}_0}{dx}=0\\
p_2\displaystyle\frac{d{a}_S}{dx}+q_2a_{S-1}=(\alpha x+\beta )a_S
\end{array}%
\right.
\end{equation}

where
$$\textbf{a}=col (a_0,\,a_1,\,...,a_S)$$ is a vector and $A,\,B$ are matrix which we determine respectively as follow
\[
\left(
\begin{array}{lllllll}
 p_1&p_0  &0      &0  &0  &\ldots&0\\
         p_2  &p_1  &p_0   &0  &0  &\ldots&0\\
         0  &p_2  &p_1   &p_0  &0  &\ldots&0\\
         \vdots&\vdots&\ldots&\vdots&\vdots&\ldots&\vdots\\
         0  &0  &0   &0  &p_2  &p_1&p_0\\
         0  &0  &0   &0  &0
         p_2&p_1,\end{array}
         \right)
         \]\\
         \[
\left(
\begin{array}{llllll}
\alpha_0-Sq_0&0  &0      &0  &0  &\ldots\\
         \alpha x+\beta-Sq_1&\alpha_0-(S-1)q_0   &0 &0  &0  &\ldots\\
         -Sq_2 &\alpha x+\beta-(S-1)q_1& \alpha_0-(S-2)q_0  &0 &0  &\ldots\\
         \vdots&\vdots&\ldots&\vdots&\vdots&\ldots\\
         0  &0    &\ldots &-2q_2 &\alpha x+\beta-q_1&\alpha_0\\
         0  &0  &\ldots   &0 &-q_2 &\alpha
         x+\beta
         \end{array}
         \right)
         \]

 From \req{57} we easily deduce that if $p_0\ne{0}$ then the
 coefficients $a_j,\quad j=0,1,2,..,S$ are polynomials of degree
 at most  $j.$

For the second case, after integration we easily deduce that
\begin{equation}\label{58}
\textbf{a}=\textsc{R}\textbf{p}
\end{equation}
 where
$$\textbf{p}=col(p^m_1,\,p^{m+1}_1,.....,p^{m+k}_1,.....,p^{m+kS}_1)$$
is a vector and $\textsc{R}$ is the following constant matrix

\[
\left(
\begin{array}{lllllllllllllll}
 R^0_0 &0    &0    &0&0     &0   &\ldots        &0        &0     &0      &\ldots    &0&\ldots   &0 &0   \\
R^1_0  &R^1_1&0    &0&0     &0     &\ldots       &R^1_k    &0        &\ldots&0      &\ldots    &\ldots  &0 &0  \\
R^2_0  &R^2_1&R^2_2&0 &0     & 0    &\ldots       &R^2_k    &R^2_{k+1}&0 &\ldots     &0&R^2_{2k}     &0  &0  \\
R^3_0  &R^3_1&R^3_2&R^3_3   &0       &  0    &\ldots&R^3_k
R^3_{k+1}&R^3_{k+2}    & 0
 0& R^3_{2k} &R^3_{2k+1}&0 \\
\vdots&\vdots&\vdots&\vdots&\vdots&0&\vdots&\vdots&\vdots&\vdots&\ldots&0
&\ldots&\ldots&
0\\
R^S_0  &R^S_1&R^S_2&\ldots&\ldots
&R^S_S&0&R^S_{k}&R^S_{k+1}&R^S_{k+2}&\ldots&0 &R^S_{2k}& \ldots&
R^S_{Sk}
\end{array}
         \right)
         \]

 Hence we easily obtain the veracity of our assertion.

\begin{cor}

{\it If the invariant algebraic curve is non reducible and
$k=\displaystyle\frac{q_0}{p_{11}}\geq{0}$ then $m=0.$}
\end{cor}
In fact, if $m\ne{0}$ and under the indicated condition we have
that the given invariant curve is reducible.

\begin{prop} {\it The maximum degree of the irreducible invariant
algebraic curve  of the non Darboux integrable quadratic system
\req{51} is $12.$ if $p_{11}\ne{0}$}
\end{prop}
In fact from \req{58},\,\req{55} by considering that the given
curve is irreducible, then we obtain that $m=0.$

On the other hand from the last of  equation of \req{57} in
particular we deduce that if $k>3$ then
\begin{equation}\label{59}
\displaystyle\frac{1}{p^2_{11}}R^S_{Sk}(p^2_{10}p_{22}-p_{10}p_{21}p_{11}+p_{20}p^2_{11})=0.
\end{equation}

Hence, if
 $R^S_{Sk}=0$ then the degree
of the algebraic curve is not maximal. On the other hand
if$$p^2_{10}p_{22}-p_{10}p_{21}p_{11}+p_{20}p^2_{11}=0$$ then the
quadratic system is Darboux integrable.

\begin{cor}

{\it Let us suppose that the algebraic curve
\begin{equation}\label{510}
g(x,y)=\sum_{l=0}^S\sum_{j=0}^{kS}\textsc{R}^l_j(p_{11}
x+p_{10})^{j+m}y^{l}=0
\end{equation}

 is invariant curve of the maximum degree
of the quadratic vector field. Then:
$$\textbf{v}=((p_{11}x+p_{10})y+p_{22}x^2+p_{21}x+p_{20})\partial_x+
(3p_{11}y^2+(q_{11}x+q_{10})y+q_{22}x^2+q_{21}x+q_{20})\partial_y$$
and}
$$dg(\textbf{v})=\displaystyle\frac{12}{p_{11}}(p^2_{11}y+p_{22}x+p_{21}p_{11}-2p_{22}p_{10})g.$$

\end{cor}
 This results can be extended analogously for the polynomial system
of degree $n.$

\begin{prop}

 {\it Let us give the invariant curve with $S$ branches of
non-Darboux integrable polynomial system  of degree $n:$
\begin{equation}\label{511}
\left\{%
\begin{array}{cc}
 \dot{x}=\sum_{j=0}^{n}p_j(x)y^{n-j}=P(x,y)\\
\dot{y}= \sum_{j=0}^{n}q_j(x)y^{n-j}=Q(x,y)
\end{array}%
\right.
\end{equation}

 with cofactor
 $K=\sum_{j=0}^{n-1}\alpha_j (x)y^{n-1-j}$
where $\alpha_j(x)=\sum_{k=0}^{j}\tau_{jk}x^k,$ and
$\tau=(\tau_{jk})$ is a constant real matrix.

Then

\begin{equation}\label{511}
\max{deg(g)}\leq \left\{%
\begin{array}{ll}
S,& \hbox{if $p_0\ne{0} $;} \\
S(n+1),& \hbox{if  $p_0={0},\, p_{11}\ne{0}$}
\end{array}%
\right.
\end{equation}

For the case when $p_{11}=0$ it is easy to  show that}
\[
\left\{%
\begin{array}{cc}
 &\alpha_0=Sq_0,\\
 &\max{deg\,a_j}\leq{j}\,,\\
 &\max{deg\,g}\leq
 {S}.
\end{array}%
\right.
\]

\end{prop}
From this result and proposition 4.2\, we deduce the proof of the
following proposition.

\begin{prop}
\, {\it The maximum degree of the invariant curve of the non
Darboux integrable polynomial planar vector field of degree $n$
($\textsc{N}(n)$) is $2n(n+1)$},\, i.e.,

\begin{equation}\label{512}
\textsc{N}(n)\leq{2n(n+1)}.
\end{equation}
\end{prop}
\section{ Quadratic system with one invariant algebraic curve.
Examples}

 First we study the case when $S=2,$ i.e., we analyze
the quadratic system with invariant algebraic curve of the type
 \begin{equation}\label{513}
g(x,y)=a_0(x)y^2+a_1(x)y+a_2(x)=0
\end{equation} where
$a_0,\,a_1,\,a_2$ are polynomial on the variable $x.$

We shall study the cases when $p_0=0,\quad p_{11}\ne{0}$ and
$p_0\ne{0}.$

For the first case the he equation \req{510}  takes the form

\begin{equation}\label{61}
g(x,y)=\sum_{l=0}^2\sum_{j=0}^{2k}\textsc{R}^l_j(p_{11}
x+p_{10})^{j}y^{l}=0.
\end{equation}
 Clearly, the maximum degree of
this curve is six.

 We shall illustrate this particular case in concrete examples.
\vskip0.25cm
 {\textbf Example 6.1.} The quadratic vector field
\[
\left\{%
\begin{array}{cc}
\dot{x}=ax^2+(-2ac+y+b)x+25C_1a^2-10aC_1d+C_1d^2\\
\dot{y}=3y^2+(-12ac+6b+dx)y+(-150C_1a^2C_3-3a^2+60C_1adC_3+\\
da-6C_1C_3d^2)x^2+(db-2adc)x- 750a^3C^2_1C_33+225C_1a^3+
450a^2C^2_1dC_3\\
-2d^3C_1-140C_1a^2d+12a^2c^2+3b^2-90d^2C^2_1aC_3-12cab+
29C_1d^2a+6d^3C_3C^2_1
\end{array}%
\right.
\]
 where $C_1,C_3$ are
arbitrary nonzero constants, admits the invariant curve of degree
six
\[
\left\{%
\begin{array}{cc}
g=C_1y^2+(x^3+(C_1d-3C_1a)x-4C_1ac+2C_1b)y+x^6C_3+(\displaystyle\frac{1}{2}d-\displaystyle\frac{3}{2}a)x^4+\\
(b-2ac)x^3+ (-3C_3C^2_1d^2-9aC_1d+
21C_1a^2-75C_3C^2_1a^2+C_1d^2+\\30C_3C^2_1ad)x^2+
(-3C_1ab-2C_1dac+C_1db+6a^2C_1c)x+\\150C_33C^3_1a^2d-\displaystyle\frac{75}{2}a^2C^2_1d+
\displaystyle\frac{15}{2}d^2C^2_1a+\\C_11b^2-250C_3C^3_1a^3-\displaystyle\frac{1}{2}d^3C^2_1+2C_3C^3_1d^3+\\
\displaystyle\frac{125}{2}a^3C^2_1-4C_1abc-30C_3C^3_1d^2a+4a^2c^2C_1=0
\end{array}%
\right.
\]

with cofactor $6(y+ax-2ac+b).$

 {\textbf Example 6.2}\,

 Now we study the case in which
$S=2,\quad p_0\ne{0}.$ Clearly, in this case the invariant
algebraic curve of the quadratic system is the family of the
conics

\begin{equation}\label{62}
  g=a_0y^2+(a_{11}x+a_{10})y+a_{22}x^2+a_{21}x+a_{20}=0
\end{equation}
 where
$a_0,\,a_{11},\,a_{10},\,a_{22},\,a_{21},\,a_{20}$ are real
constants.

In particular, for the quadratic system
 \begin{equation}\label{63}
\left\{%
\begin{array}{cc}
 \dot{x}=\beta\,y^2-(\beta+2)xy+(\beta-4)x^2-2\beta x\\
\dot{y}=(\alpha-4)y^2-2(\alpha+2)xy+(\alpha-4)x^2-2\alpha\,x
\end{array}%
\right.
\end{equation}
 we have that $p_0=\beta\ne{0}.$

 After some calculations we can prove that the invariant curve is the
 parabola
$$(y-x)^2-2x=0$$
with the cofactor
$$K=2(\beta+\alpha)y+2(\beta+\alpha)x-2\beta .$$

The quadratic system with invariant parabola was constructed  in
[Sad1] and admits the equivalent representation :
\begin{equation}\label{64}
\left\{%
\begin{array}{cc}
 \dot{x}=\beta\,\big((y-x)^2-2\,x\big)-4\,x\,(x + y)\\
 \dot{y}=\alpha\,\big((y-x)^2-2\,x\big)-4\,(x + y)^2
\end{array}%
\right.
\end{equation}
The points
$$O(0,\,0),\quad N(\displaystyle\frac{1}{2},\,-\displaystyle\frac{1}{2}),\quad M(\displaystyle\frac{\beta^2}{K_3},\,
\displaystyle\frac{\beta^2(2\alpha-\beta)}{K_3})$$ are its
critical points, where $K_3\equiv{2((\alpha-\beta)^2-2\alpha}).$
 The bifurcation analysis show that this system is generic.
 The bifurcation curves divide the plane $(\alpha,\,\beta)$ in 17
 region in which we observe a qualitative change in the behavior
 of the trajectories of the constructed quadratic system. We
 determine 38 different quadratic systems, among these there is
 one with one limit cycles\,\cite{Rey}.

The quadratic system
\[
\left\{%
\begin{array}{cc}
\dot{x}=\displaystyle\frac{1}{2A}(Ax^2+By^2+2C)A_0+\displaystyle\frac{1}{A}(-2Byq_{11}-2xBq_{22}+xA_{1}A_2Bq_{21})y\\
\dot{y}=\displaystyle\frac{1}{2B}(A_1y^2+yq_{11}x+q_{22}x^2+q_{21}x+A_1C),
\end{array}%
\right.
\]
admits as invariant the curve
$$Ax^2+By^2+2C=0$$ with cofactor
$K=A_{0}x+A_{1}y,$ where $A,B,C,A_0,A_1,q_{11},q_{21},q_{22},$ are
real constants such that $A\ne{0},\,B\ne{0}.$

 Now we study the case when $S=3,$ i.e., we analyze the
 quadratic system with invariant algebraic curve of the type

\begin{equation}\label{65}
 g(x,y)=\sum_{l=0}^2\sum_{j=0}^{3k}\textsc{R}^l_j(p_{11} x+p_{10})^{j}y^{l}=0.
\end{equation}

 {\textbf Example 6.3}\quad (The Filipstov system ).

For the quadratic system

\[
\left\{%
\begin{array}{cc}
\dot{x}=16(1+a)x-6(2+a)x^2+(2+12x)y\\
\dot{y}=3a(1+a)x^2+(15(1+a)-2(9+5a)x)y+16y^2
\end{array}%
\right.
\]

 we have that
$3k=4\Longrightarrow{m=0}.$

  As we can observe in this case the quadratic system
possesses the irreducible invariant algebraic curve
$$g(x,y)\equiv{y^3+\displaystyle\frac{1}{4}(3(1+a)-6(1+a)x)y^2+\displaystyle\frac{3}{4}(1+a)ax^2y+\displaystyle\frac{3}{4}(1+a)a^2x^4}=0.$$
The cofactor of this curve is $$K=48y-4(1+a)x+5(1+a).$$

{\textbf Example 6.4}\quad For the quadratic system
\[
\left\{%
\begin{array}{cc}
\dot{x}=(2+3a)x^2+(2+4y)x+y\\
\dot{y}=(5+4(1+a))y+6y^2+ax^2
\end{array}%
\right.
\]
 we obtain that
$k=2\Longrightarrow{m=0.}$

The quadratic system possesses the irreducible invariant algebraic
curve of degree four \cite{Chr et al.3}
$$a^2x^4++2ax^2(x+1)y+(1+x)y^2+y^3=0$$

The cofactor $K$ in this case is
$$K=18y+(5+6a)x+5.$$

 For the case when $S=4$ we obtain the invariant algebraic curve of
the type

\begin{equation}\label{66}
   g(x,y)=\sum_{l=0}^2\sum_{j=0}^{4k}\textsc{R}^l_j(p_{11} x+p_{10})^{j}y^{l}=0.
 \end{equation}
Clearly, the maximum degree of this curve is 12. The upper bound
is reached in particular in the following example  \cite{Chr et
al.3}

{\textbf Example 6.5}\, For the non Darboux integrable quadratic
system
\[
\left\{%
\begin{array}{cc}
\dot{x}=xy+x^2+1\\
\dot{y}=3y^2-\displaystyle\frac{81}{2}x^2+\displaystyle\frac{57}{2}
\end{array}%
\right.
\]

has we have that $k=3\Rightarrow{m=0}$ as a consequence the above
vector field has the invariant irreducible algebraic curve of the
maximum  degree 12.  In the indicated paper was showed that the
the curve
\[
\left\{%
\begin{array}{cc}
-442368-7246584x^2+71546517x^4-97906500x^6+41343750x^8-23437500x^{10}+\\
48828125x^{12}+(322272x-12126312x^3+23463000x^5+1125000x^7+15625000x^9)y-
\\(98784-711288x^2+5058000x^4-375000x^6)y^2+(32928x-1124000x^3)y^3-5488y^4=0
\end{array}%
\right.
\]

is its invariant.

\section{ Construction the polynomial planar system
 with invariant algebraic curves with variables
separable}

In this section we deal with the polynomial system with invariant
algebraic curve with variables separable
$$g(x,y)=F_1(x)+F_2(y)=0,$$ where $F_1,\,F_2$ are arbitrary
polynomials :
$$deg(g(x,y))=max\Big(deg(F_1(x)),\,max(deg(F_2(y))\Big).$$

We state and study the following problem.

 {\textbf Problem 7.1 }\,Let $g$ be  a function:
$$g(x,y)=g_0+A\int\prod_{j=1}^{m_1} (x-a_j)dx+B\int\prod_{j=1}^{m_2}
(y-b_j)dy$$ where
$A,\,B,\,a_1,\,a_2,\,...,a_{m_1},\,b_1,\,b_2,\,...,\,b_{m_2}$ are
real parameters such that
$$a_1<\,a_2<....<a_{m_1},\quad b_1<b_2.....<b_{m_2},\quad AB\ne{0}.$$
We require to determine the non-Darboux integrable  polynomial
vector  $\textbf{v}$ of degree $n$ for which the given curve is
its invariant.

We propose the solution of the state problem  for the following
particular cases
\[
\left\{%
\begin{array}{cc}
m_2=n-1,\,m_1=n-1,\\
\,m_2=2m-1,\,m_1=2m+1,\,n=2m+1,\\
 m_2=2m-2,\,m_1=2m,\,n=2m,\\
m_2=m_1=m,\quad n=2m+1, \\
m_2=m_1=m,\quad n=2m+2.
\end{array}%
\right.
\]

\begin{prop}
 \ The polynomial system of degree $n$ \cite{Sad2}
\[
\left\{%
\begin{array}{cc}
\dot{x}=(Ax+By+C)\partial_yg(x,y)\equiv{P(x,y)}\\
\dot{y}=-(Ax+By+C)\partial_xg(x,y)+\lambda g(x,y)\equiv{Q(x,y)}
\end{array}%
\right.
\]
admits as invariant curve
$$g(x,y)=g_0+K_1\int\prod_{j=1}^{n-1}(y-b_j)dy+ K_2\int \prod_{j=1}^{n-1}(x-a_j)dx=0,$$
\end{prop}
where
$a_1,\,a_2,\,...,a_{n-1},\,b_1,\,b_2,...,b_{n-1},\,K_1,\,K_2,\,g_0,\,A,\,B,\,C,\,\lambda$
are  arbitrary real parameters.

By choosing the arbitrary parameters properly we can construct the
nonsingular algebraic curve of degree $n,$ hence the genus
($\textsc{G}$) of this curve is:
$$\textsc{G}=\displaystyle\frac{1}{2}(n-1)(n-2)$$

{\textbf Example 7.1}\quad Let $g$ is a nonsingular curve of
degree $2m+2$ such that
$$g_m(x,y)=g_o+\int\prod_{j=1}^{m}x((\displaystyle\frac{x}{j\pi})^2-1)dx+\int\prod_{j=1}^{m}y((\displaystyle\frac{y}{j\pi})^2-1)dy=0$$
It is easy to show that the polynomial system of degree $n=2m+1:$

\[
\left\{%
\begin{array}{cc}
\dot{x}=(A_mx+B_my+C_m)\prod_{j=1}^{m}y((\displaystyle\frac{y}{j\pi})^2-1)\equiv{P(x,y)}\\
\dot{y}=-(A_mx+B_my+C_m)\prod_{j=1}^{m}x((\displaystyle\frac{x}{j\pi})^2-1)+\lambda
g_m(x,y)\equiv{Q(x,y)}
\end{array}%
\right.
\]

admits as invariant the given curve.

By considering that
$$\prod_{j=1}^{\infty}((\displaystyle\frac{y}{j\pi})^2-1)=\sin y\Rightarrow\lim_{m\to{+\infty}}g_m(x,y)=g_0+\cos x+\cos y$$
and choose the arbitrary parameters  $A_m,\,B_m,\,C_m$ properly,
we obtain the analytic planar vector field

\[
\left\{%
\begin{array}{cc}
\dot{x}=\sin y\equiv{P(x,y)}\\
\dot{y}=-\sin x+\lambda (g_0+\cos x+\cos y)\equiv{Q(x,y)}
\end{array}%
\right.
\]

for which the curve
$$g_0+\cos x+\cos\,y=0$$ is its invariant.

Clearly, the constructed analytic system admits infinity many
number of limit cycles. \vskip0.25cm
 \begin{prop}

 \quad The  polynomial vector field of degree $n$

\begin{equation}\label{71}
\left\{%
\begin{array}{cc}
\dot{x}=(a+byx)\partial_yH(x,y)\\
\dot{y}=-(a+byx)\partial_xH(x,y)+(n+1)by H(x,y)
\end{array}%
\right.
 \end{equation}

 admits as
invariant the algebraic curve of degree $n+1$

\begin{equation}\label{72}
  H(x,y)\equiv{x^{n+1}+G_{n-1}(x,y)}=0,
\end{equation}
 where $G_{n-1}$ is an arbitrary polynomial of degree
$n-1.$
\end{prop}

Clearly, this system in general  has no Darboux integrating
factors or first integrals.

 The following particular case is an interesting one:

\begin{cor}

 Let
$$g_m(x,y)=g_0+\int_{x_0}^x\prod_{j=1}^{m+1}x(x^2-a^2_j)dx+\int_{y_0}^y\prod_{j=1}^{m-1}y(y^2-b^2_j)dy=0$$
is a curve of degree $2m+2$ with the maximum genus
$\textsc{G}={2(m+1)(m-1)+1}$ is invariant of the vector field of
degree $n=2m+2:$
\[
\left\{%
\begin{array}{cc}
\dot{x}=(a+b_myx)\prod_{j=1}^{m-1}y(y^2-b^2_j)\\
\dot{y}=-(a+b_myx)\prod_{j=1}^{m+1}x(x^2-a^2_j+(2m+2)b_my g_m(x,y)
\end{array}%
\right.
\]

\end{cor}

{\textbf Example 7.2} By making $m\to{+\infty}$ and choose the
arbitrary parameters properly we deduce from the above system as a
particular case the analytic system
\[
\left\{%
\begin{array}{cc}
\dot{x}=aJ_0(y)\\
\dot{y}=-aJ_0(x)+\lambda y (J_1(x)+J_1(y)+g_0)
\end{array}%
\right.
\]

where $J_0,\,J_1$ are the Bessel functions. This analytic system
admits an infinity many number of limit cycles. \vskip0.25cm
 Analogously we construct the
polynomial system of degree $n=2m:$
 \[
\left\{%
\begin{array}{cc}
\dot{x}=&(a+b_myx)\prod_{j=1}^{m-1}(y^2-b^2_j)\\
\dot{y}=&-(a+b_myx)\prod_{j=1}^{m+1}(x^2-a^2_j)+2mb_my g_m(x,y)
\end{array}%
\right.
\]

with invariant curve
$$g_m(x,y)=g_0+\int_{x_0}^x\prod_{j=1}^{m+1}(x^2-a^2_j)dx+\int_{y_0}^y\prod_{j=1}^{m-1}(y^2-b^2_j)dy=0$$
is a curve of degree $2m+2$ with maximum genus
$\textsc{G}={2m(m-1)}$

By using the algebraic packages it is possible to show the
following proposition.

\begin{prop}

 There exist polynomials $p(x,y),\,q(x,y)$ of degree $n$ for
which the non-Darboux integrable differential system
\begin{equation}\label{73}
\left\{%
\begin{array}{cc}
\dot{x}=a\partial_yg(x,y)+p(x,y),\quad a=const.\\
\dot{y}=-a\partial_xg(x,y) +q(x,y)
\end{array}%
\right.
 \end{equation}

has the invariant curve of degree $n+1$
\begin{equation}\label{74}
g(x,y)=g_0+A\int(\prod_{j=1}^{n}(x-a_j))dx+B\int(\prod_{j=1}^{n}
(y-b_j))dy=0,
\end{equation}
 for certain values of the real parameters
$g_0,A,B,\,a_1,...,a_m,\,b_1,b_2,...,b_m,a_0,\,b_0.$
\end{prop}
Clearly if this curve is non singular then the genus is
$\textsc{G}=\displaystyle\frac{1}{2}n(n-1).$ \vskip0.25cm
 {\textbf Example 7.3}

Let us suppose that the given algebraic curve (7.4) is such that
$$g(x,y)=g_0+A\int(\prod_{j=1}^{m}x(x^2-a^2_j))dx+B\int(\prod_{j=1}^{m}
y(y^2-b^2_j))dy.$$

It is possible to construct the non-Darboux integrable polynomial
vector fields of degree $n=2m+1.$ In particular for $n=3,5,7$ we
construct the following non-Darboux integrable polynomial systems.

For the polynomial system of degree seven
\[
\left\{%
\begin{array}{cc}
\dot{x}=
y(160p^4y^2q^2-192p^2y^4q^2+64p^2y^6-96p^4y^4+32p^6y^2-64p^4q^4-32p^6q^2-\\
96y^4q^4-32p^2q^6+32y^2q^6+160p^2y^2q^4+64y^6q^2)\nu_0+\lambda y(-12p^2x^4-4p^4y^2-4y^2q^4+\\
8p^2y^2x^2+4p^2y^4+q^6-p^2q^4+12x^2p^2q^2+8y^2x^2q^2-8p^2y^2q^2+8x^6+2x^2p^4+\\
2x^2q^4-8y^4x^2-12x^4q^2-p^4q^2+p^6+4y^4q^2)\\
\dot{y}=-1/64x(-64p^4q^4-32p^2q^6-96x^4p^4+32p^6x^2+32x^2q^6-32p^6q^2-\\
192x^4p^2q^2+160p^4x^2q^2+160p^2x^2q^4+64p^2x^6-96x^4q^4+64x^6q^2)/(p^2+q^2)\nu_0\\-1/64x\lambda
(8p^2y^2x^2-12p^2y^4-8y^2x^4-p^2q^4-8x^2p^2q^2+8y^2x^2q^2+
12p^2y^2q^2+2p^4y^2-\\4x^2p^4-4x^2q^4+q^6+
p^6+2y^2q^4-p^4q^2+4p^2x^4+4x^4q^2- 12y^4q^2+8y^6)/(p^2+q^2)
\end{array}%
\right.
\]
 the invariant curve  is
\[
\left\{%
\begin{array}{cc}
g(x,y)=1/8x^8+(-1/4p^2-1/4q^2)x^6+(1/8p^4+1/2p^2q^2+1/8q^4)x^4-\displaystyle\frac{1}{4}(p^4q^2+p^2q^4)x^2+\\1/8y^8+
(-1/4p^2-1/4q^2)y^6+(1/8p^4+1/2p^2q^2+1/8q^4)y^4+(-1/4p^4q^2-\\1/4p^2q^4)y^2-1/128q^8+1/32q^6p^2+13/64q^4p^4+1/32q^2p^6-1/128p^8=0
\end{array}%
\right.
\]

with the cofactor
$$K=\lambda (-y+x)(x+y)(-p^2-q^2+y^2+x^2)yx.$$
In this example we have that $$a_1=-a_2=p,\quad a_3=-a_4=q,\quad
a_5=-a_6=\sqrt{\displaystyle\frac{1}{2}(p^2+q^2)}.$$  For the
quintic vector field
\[
\left\{%
\begin{array}{cc}
\dot{y}=\nu_0(-3x^2+r^2)(6r^2-6x^2)x+\lambda (-3x^2+r^2)(-r^2+y^2)y\\
\dot{x}=-\nu_0(6r^2-6y^2)(r-3y^2)y-\lambda (r^2-x^2)x(r^2-3y^2)
\end{array}%
\right.
\]

The curve
$$g(x,y)=\displaystyle\frac{1}{6}y^6-\displaystyle\frac{1}{3}r^2y^4+\displaystyle\frac{1}{6}r^4y^2+\displaystyle\frac{1}{6}x^6-\displaystyle\frac{1}{3}x^4r^2+
\displaystyle\frac{1}{6}r^4x^2-\displaystyle\frac{1}{39} r^6=0$$
is its invariant  with the cofactor
$$K=6y^2\lambda (3x^2-r^2).$$

The cubic polynomial system

\[
\left\{%
\begin{array}{cc}
 \dot{x}=\nu_0( y^3-y)+\lambda x(2q^2x^2-2p^2y^2-3q^4+p^4)\\
\dot{y}=-\nu_0(x^3-x)+\lambda
y(2q^2x^2-2p^2y^2+3q^4-p^4),
\end{array}%
\right.
\]
 has as invariant
curve of degree four
$$g(x,y)= 1/8q^4+1/8p^4+1/4x^4-1/2q^2x^2+1/4y^4-1/2p^2y^2=0$$
with cofactor
$$ K=4\lambda (2q^2x^2-2p^2y^2-q^4+p^4).$$

It is interesting to observe that the above constructed polynomial
system of degree seven under the change
\[
\left\{%
\begin{array}{cc}
 x=\sqrt{X},\quad y=\sqrt{Y},\\
\bold{V}=\Big(\displaystyle\frac{\textbf{v}(x)}{y},\,\displaystyle\frac{\textbf{v}(y)}{x}\Big)|_{x=\sqrt{X},y=\sqrt{Y}},
\end{array}%
\right.
\]
can be transformed to the cubic system
\[
\left\{%
\begin{array}{cc}
\dot{Y}=(-64p^2-64q^2)X^3+(96q^4+192p^2q^2+96p^4)X^2+\\
(-32p^6-32q^6-160q^2p^4-160q^4p^2)X+64p^4q^4+32p^2q^6+32p^6q^2)p0+(-8Y^3+\\
(12p^2+12q^2)Y^2+(-2p^4-2q^4-8p^2X-12p^2q^2-8Xq^2+8X^2)Y+q^4p^2-p^6-q^6+\\
q^2p^4+8Xp^2q^2-4p^2X^2+4Xq^4+4Xp^4-4X^2q^2
)\lambda)\\
\dot{X}=-((-64p^2-64q^2)Y^3+(96q^4+192p^2q^2+96p^4)Y^2+\\
(-32p^6-32q^6-160q^2p^4-160q^4p^2)Y+64p^4q^4+32p^2q^6+32p^6q^2)p0+8X^3+\\
(12p^2+12q^2)X^2+(-12p^2q^2-8p^2Y-8Yq^2+8Y^2-2q^4-2p^4)X-\\
p^6-q^6+q^4p^2+q^2p^4-4Y^2q^2+4Yq^4+8p^2Yq^2+4p^4Y-4p^2Y^2)\lambda
)\end{array}%
\right.
\]

 which admits the invariant curve of degree four
\[
\left\{%
\begin{array}{cc}
g(X,Y)=1/8X^4+(-1/4p^2-1/4q^2)X^3+(1/2p^2q^2+1/8p^4+1/8q^4)X^2+\\
(-1/4q^2p^4-1/4q^4p^2)X+1/8Y^4+(-1/4p^2-1/4q^2)Y^3+(1/2p^2q^2+\\
1/8p^4+1/8q^4)Y^2+(-1/4q^2p^4-1/4q^4p^2)Y-1/128p^8+1/32p^6q^2-\\
1/128q^8+13/64p^4q^4+1/32p^2q^6
\end{array}%
\right.
\]

 with cofactor
$$K(X,Y)=(X-Y)(X+Y-(p^2+q^2))\lambda .$$

{\textbf Example 7.4}

 Finally we analyze the case when
$$g(x,y)=g(x_0,y_0)+\int_{x_0}^x\prod_{j=1}^{m+1}x(x^2-a^2_j)dx+\int_{y_0}^y\prod_{j=1}^{m}y(y^2-b^2_j)dy$$
is a curve of degree $2m+4.$

 Analogously to the above case we can construct a non-Darboux
integrable polynomial vector fields of degree $n=2m+2$

In particular for $n=4,6$ we construct the following polynomial
systems.

 For $n=4$
 \[
\left\{%
\begin{array}{cc}
\dot{x}=(1/6y^2x^2-5/108y^2-1/18x^4+7/324x^2+1/729)\nu_0+(1/6xy^3-1/54xy)\lambda\\
\dot{y}=(-1/54xy+1/4xy^3-1/12yx^3)\nu_0+(1/4y^4-5/108x^4-1/18y^2+4/2187+4/243x^2)\lambda
\end{array}%
\right.
\]

the invariant curve is
$$g(x,y)=
1/6x^6-5/36x^4+2/81x^2+1/4y^4-1/18y^2+4/2187=0.$$   The cofactor
in this case
is$$K=-\displaystyle\frac{1}{236196}(1458x^6-1215x^4+216x^2+2187y^4-486y^2+16)\lambda
.$$

For $n=6$ we construct the following vector field
\[
\left\{%
\begin{array}{cc}
\dot{x}=(-1/128y^4x^2-4/3y^2+1/64x^6-1/2x^4+1/12y^4+40/9x^2+1/8y^2x^2-\\256/27)\nu_0+
&(1/8xy^3-3/512xy^5-1/2xy)\lambda\\
\dot{y}=(-\displaystyle\frac{1}{96}xy^5+\displaystyle\frac{5}{18}xy^3+\displaystyle\frac{1}{48}yx^5-\\
\displaystyle\frac{4}{9}yx^3+\displaystyle\frac{16}{27}xy)\nu_0+
(-\displaystyle\frac{1}{128}y^6\\-\displaystyle\frac{11}{6}x^4+\displaystyle\frac{1}{4}y^4-\displaystyle\frac{1024}{27}-2y^2+16x^2+
\displaystyle\frac{1}{16}x^6)\lambda
\end{array}%
\right.
\]
The invariant curve and its cofactor are respectively
$$g(x,y)=1/6y^6-16/3y^4+128/3y^2+1/8x^8-16/3x^6+704/9x^4-4096/9\,x^2+65536/81=0,$$

$$ K=(-4y+y^3-3/64y^5)\lambda+(1/8\,x^5-8/3\,x^3-1/16\,y^4\,x+y^2x+32/3\,x)\nu_0 .$$

We observe that after the change
\[
\left\{%
\begin{array}{cc}
 x={X}^2,\quad y={Y}^2,\\
\bold{V}=\Big(\displaystyle\frac{\textbf{v}(x)}{x},\,\displaystyle\frac{\textbf{v}(x)}{y}\Big)|_{x={X}^2,\,y={Y}^2}
\end{array}%
\right.
\]

the above system takes the form
\[
\left\{%
\begin{array}{cc}
\dot{X}=Y((-1/2X^2Y^2+1/8X^2Y^6-3/512X^2Y^{10})\nu_0+\\
(-1/128Y^8X^4+1/64X^{12}-256/27-4/3Y^4+40/9X^4+1/8Y^4X^4-1/2X^8+1/12Y^8)\lambda) \\
\dot{Y}=X((-1/128Y^12-2Y^4+1/4Y^8+1/16X^{12}-11/6X^8+16X^4-1024/27)\nu_0+\\
(-1/96X^2Y^{10}-4/9Y^2X^6+5/18X^2Y^6+1/48Y^2X^{10}+16/27X^2Y^2)\lambda)
\end{array}%
\right.
\]

This polynomial system of degree 13 admit as invariant curve of
degree 16

$$g(X,Y)=1/6Y^{12}-16/3Y^8+128/3Y^4+1/8X^{16}-16/3X^{12}+704/9X^8-4096/9X^4+65536/81=0$$

with cofactor

$$K=\displaystyle\frac{1}{96}(\nu_0(24X^{10}-512X^6+192Y^4X^2-12Y^8X^2+2048X^2)+\lambda(192Y^6-9Y^{10}-768Y^2)YX)$$

We observe that the constructed above curve are singular curves.

We state the following problem

{\textbf Problema }\quad To construct the differential system
(7.3) for which the invariant curve (7.4) is non singular.

Clearly, if the curve (7.4) is non singular then the genus is
$$\textsc{G}=\displaystyle\frac{1}{2}n(n-1).$$
as a consequence the maximum number of algebraic limit cycles are
$\displaystyle\frac{1}{2}n(n-1)+1.$

To construct the require vector field should be satisfied  the
relation $dg(\textbf{v})=K(x)g$ under the condition that the curve
$g(x,y)=0$ is non singular. To obtain the explicit expression for
the vector field in general it is necessary to solve a lot of
technical and theoretical problems. In particular for the cubic
system we have 31 parameters which must be satisfies 27 equations.
By solving these equations we obtain that the nonsingular curve of
genus 3:
\[
\left\{%
\begin{array}{cc}
g(x,y)=1/4x^4-1/3(a_2+a_1)x^3+1/2a_2a_1x^2+1/4y^4-1/3(b_2+b_1)y^3+1/2b_2b_1y^2+\\1/12\mu\Big(-5b_2^3b_1^3a_2^2-
5b_2^3b_1^3a_1^2+2b_2^4b_1^2a_1^2+2b_2^2b_1^4a_1^2-2b_1a_1^2a_2^4b_2+2b_2^4b_1^2a_2^2+\\2b_2^2b_1^4a_2^2-
2b_2^2b_1^4a_1a_2+2a_1^4a_2^2b_2^2+5b_2^3b_1^3a_1a_2-
2b_2^4b_1^2a_1a_2+\\5b_1a_1^3a_2^3b_2+2b_1^2a_1^2a_2^4-5a_1^3a_2^3b_2^2-2b_1a_1^4a_2^2b_2+2a_1^2a_2^4b_2^2-5b_1^2a_1^3a_2^3+
2b_1^2a_1^4a_2^2\Big)=0,\\
\mu\equiv{\Big((a_1^2-a_1a_2+a_2^2)(-b_2b_1+b_2^2+b_1^2)\Big)^{-1}}
\end{array}%
\right.
\]
which is invariant of the cubic system if the parameters
$a_1,\,a_2,\,b_1,\,b_2$ are roots of the certain homogenous
polynomials $P_j(a_1,a_2,b_1,b_2),\ j=1,2,3$ such that
 \[
\left\{%
\begin{array}{cc}
P_j(\lambda
a_1,\lambda a_2,\lambda b_1,\lambda
b_2)=\lambda^8P_j(a_1,a_2,b_1,b_2),\quad j=1,2\\
P_3(\lambda a_1,\lambda a_2,\lambda b_1,\lambda
b_2)=\lambda^{24}P_3(a_1,a_2,b_1,b_2)
\end{array}%
\right.
\]

\section{The Poincar\'e problem and 16th Hilbert problem  for algebraic limit cycles}
 {\textbf 8.1\quad The Poincar\'e problem.}

The question on the existence an effective procedure to find a
natural number $\textsc{N}(n)$ which bounds the degree of all
irreducible invariant curve of a non-Darboux integrable polynomial
system of a degree $n$ is well known as Poincar\'e problem
\cite{Poin},\cite{Car},\cite{Sad1}.

 The problem on the existence  of the
polynomial  planar systems with an invariant algebraic curve of
the  maximum degree was studied in
 particular in \cite{Chr},\\cite{Sad2}.  In \cite{Chr} the author  gave an explicit polynomial system of
degree $n$ for each non-singular real algebraic curve $g=0$ of
degree $n$ what is that system's invariant. The author states also
that $n$ is optimal for a generic class of algebraic curve.

\begin{prop}

Let $\textsc{N}(n)$ is the maximum degree of the irreducible
algebraic curve (4.9,\,(4.10) that is invariant curve of the
polynomial system (3.5) of degree $n$.

Then $$(n+1)\leq\textsc{N}(n)\leq 2n(n+1)$$
\end{prop}
The proof of the lower bound follows from the proposition 7.2 and
7.3 and the upper bound from the proposition 5.4. The upper bound
is reached in particular for $n=2.$

 \section{ The 16th Hilbert Problem for Algebraic Limit Cycles}

 In 1900 Hilbert\cite{Hil} proposed in the
second part of his 16th problem to estimate a uniform upper bound
for the number of limit cycles of all polynomial vector fields of
a given degree. This question has been studied in particular in
\cite{Sad1},\cite{Sad2} for algebraic limit cycles.

 By considering that the ovals of the
invariant algebraic curve are isolated periodical solutions of the
vector field for which is it invariant we deduce that the maximum
number of algebraic limit cycles of the polynomial system of
degree $n$ with one invariant curve (we denote this number  by
$A(n,1)$ ) is at most the {\it genus of the curve +1,} i.e.,
$$A(n,1)\leq{\textsc{G}+1}$$
 hence
we observe that to solve the 16th Hilbert's Problem for Algebraic
Limit Cycles it is necessary firstly to solve the Poincar\'e
problem for this case, i.e.,it is necessary to find the maximum
degree of the invariant curve \cite{Sad1}.

\begin{prop}

$$A(n,1)\leq{(2n^2+2n-1)(n^2+n-1)+1}.$$
\end{prop}
The proof follows from the proposition  (5.4) and from the Harnak
theorem on the maximum numbers of the ovals of the algebraic
curve.

From the above results we prove the following

 \begin{prop}

{\it The maximum number of the algebraic limit cycles  for the
polynomial planar vector field of degree $n$ with one invariant
algebraic curve $A(n,1)$ is such that}

$$A(n,1)\geq{max\big(h_0(n),h_1(n),h_2(n),h_3(n),h_4(n),\big)}$$
\end{prop}

where $h_j(n)$ are the maximum numbers of ovals of the algebraic
curves  $H_j(x,y)=0:$
\[
\left\{%
\begin{array}{cc}
H_0(x,y)=\sum_{m=0}^{n}a_{n-m\,m}x^my^{n-m},\\
H_1(x,y)=ax^{n+1}+\sum_{m=0}^{n-1}a_{n-m\,m}x^my^{n-m},\\
H_k(x,y)=\int(\prod_{1}^{m_k}(x-a_j)dx+\prod_{1}^{l_k}(y-b_j)dy),\\
a_1<a_2....<a_{m_
k},\quad b_1<b_2....<b_{l_k},\quad k=2,3,4\\
 m_2=l_2=n-1,\quad m_3=n+1,\quad l_3=n-1,\quad
m_4=n,\,l_4=n.
\end{array}%
\right.
\]

From this inequality we deduce the following result \cite{Sad2};

 \begin{cor}

$$A(n,1)\geq{\displaystyle\frac{1}{2}(n-1)(n-2)}+1.$$
\end{cor}

{\textbf Conjecture}

$$A(n,1)\geq{\displaystyle\frac{n(n-1)}{2}}$$
 This conjecture can be solve if we construct the polynomial
system of degree $n$ with one invariant nonsingular irreducible
invariant curve of degree $n+1$ (see proposition 8.1).

 It is interesting to observe
that in \cite{Llib1} stated the following problems:

Is 1 the maximum number of algebraic limit cycles that a quadratic
system can have?.

Is 2 the maximum number of algebraic limit cycles that a cubic
system can have?.

The answer to this questions for the system (3.5) with $S>1$ was
given in \cite{Sad2}.

Is there a uniform bound for the number of algebraic limit cycles
that a polynomial vector field of degree $n$ could have?.
\vskip0.25cm

 What is the maximum degree of an algebraic limit cycle of a
quadratic polynomial vector field?.

 What is the maximum degree of an algebraic limit cycle of a
cubic polynomial vector field?

The partial solutions of these problems we obtain from the above
results and results proposed in \cite{Sad1},\,\cite{Sad2}.

\vskip0.25cm

 {\textbf Acknowledgements}

\vskip0.25cm

 This work was partly supported by the Spanish Ministry of
Education through projects DPI2007-66556-C03-03,
TSI2007-65406-C03-01 "E-AEGIS" and Consolider CSD2007-00004
"ARES".

\end{document}